\documentclass[12pt]{article}

\usepackage{amssymb}
\usepackage[mathscr]{eucal}
\usepackage{amsmath,amssymb,latexsym,theorem,bbm}

\newfont{\msa}{msam10 scaled\magstep1}
\newfont{\ssmsa}{msam9}
\newfont{\smsa}{msam10}
\newfont{\sms}{msbm10}
\newfont{\sseufb}{eufb9}
\newfont{\seufb}{eufb10}
\newfont{\eufb}{eufb10 scaled\magstep1}
\newfont{\eusb}{eusb10 scaled\magstep1}
\newfont{\hcmr}{cmr17 scaled\magstep5}

\newcommand{\DS}{\displaystyle}

\newcommand{\CC}{\mathbb{C}}

\newcommand{\EE}{\mathsf{E}}
\newcommand{\HH}{\mathbb{H}}
\newcommand{\GG}{\mathbb{G}}

\newcommand{\NN}{\mathbb{N}}

\newcommand{\RR}{\mathbb{R}}

\newcommand{\ZZ}{\mathbb{Z}}

\newcommand{\sSS}{\raise-0.5truemm\hbox{\sms S}}

\newcommand{\cA}{{\cal A}}
\newcommand{\cG}{{\cal G}}

\newcommand{\cF}{{\cal F}}
\newcommand{\cH}{{\cal H}}
\newcommand{\cL}{{\cal L}}

\newcommand{\dd}{\mathsf{d}}

\newcommand{\ee}{\mathsf{e}}
\newcommand{\cov}{\mathsf{Cov}}
\newcommand{\var}{\mathsf{Var}}

\newcommand{\rank}{\mathsf{rank}\,}
\renewcommand{\Im}{\mathsf{Im}\,}
\renewcommand{\Re}{\mathsf{Re}\,}

\newcommand{\Det}{\mathsf{det}\,}
\newcommand{\supp}{\mathsf{supp}\,}

\newcommand{\hmu}{\widehat{\mu}}

\newcommand{\tf}{\widetilde{f}}
\newcommand{\tC}{\widetilde{C}}
\newcommand{\tF}{\widetilde{F}}

\newcommand{\ts}{\widetilde{s}}

\newcommand{\txi}{\widetilde{\xi}}
\newcommand{\teta}{\widetilde{\eta}}

\newcommand{\sleq}{\mbox{\ssmsa\hspace*{0.1mm}\symbol{54}\hspace*{0.1mm}}}
\newcommand{\sgeq}{\mbox{\ssmsa\hspace*{0.1mm}\symbol{62}\hspace*{0.1mm}}}
\renewcommand{\leq}{\mbox{\msa\hspace*{0.9mm}\symbol{54}\hspace*{0.9mm}}}
\renewcommand{\geq}{\mbox{\msa\hspace*{0.9mm}\symbol{62}\hspace*{0.9mm}}}

\newcommand{\bone}{\mathbbm{1}}
\newcommand{\vare}{\varepsilon}

\newcommand{\proofend}{\hfill\mbox{$\Box$}}

\theorembodyfont{\em}
\newtheorem{Lem}{Lemma}[section]
\newtheorem{Thm}{Theorem}[section]
\newtheorem{Pro}{Proposition}[section]

\theorembodyfont{\rm}
\newtheorem{Rem}{Remark}[section]

\numberwithin{equation}{section}

\newcommand{\sphat}%
 {\hspace*{-0.8mm}\raise1.5truemm\hbox{$\widehat{ }$}\raise-1.5truemm\hbox{}\,}
\newcommand{\ssphat}%
 {\hspace*{-0.7mm}\raise1truemm\hbox{$\widehat{ }$}\raise-1truemm\hbox{}\,}

\begin{document}

\title{\bfseries\Large Fourier transform of a Gaussian measure on the
                        Heisenberg group }

\author{ {\sc M\'aty\'as Barczy} and {\sc Gyula Pap} \\
           University of Debrecen, Hungary}

\date{}

\maketitle

\renewcommand{\thefootnote}{}

\vspace*{-10mm}

\begin{abstract}
An explicit formula is derived for the Fourier transform of a Gaussian measure
 on the Heisenberg group at the Schr\"odinger representation.
Using this explicit formula, necessary and sufficient conditions are
 given for the convolution of two Gaussian measures to be a Gaussian measure.
\end{abstract}

\section{Introduction}

Fourier transforms of probability measures on a locally compact topological
 group play an important role in several problems concerning convolution and
 weak convergence of probability measures.
Indeed, the Fourier transform of the convolution of two probability measures
 is the product of their Fourier transforms, and in case of many groups the
 continuity theorem holds, namely, weak convergence of probability measures
 is equivalent to the pointwise convergence of their Fourier transforms.
Moreover, the Fourier transform is injective, i.e., if the Fourier transforms
 of two probability measures coincide at each point then the measures coincide.
(See the properties of the Fourier transform, e.g., in Heyer
 \cite[Chapter I.]{HEY77}.)
In case of a locally compact Abelian group, an explicit formula is available
 for the Fourier transform of an arbitrary infinitely divisible probability
 measure (see Parthasarathy \cite{PAR67}).
The case of non-Abelian groups is much more complicated.
For Lie groups, Tom\'e \cite{TOM98} proposed a method how to calculate Fourier
 transforms based on Feynman's path integral and discussed the physical
 motivation, but explicit expressions have been derived only in very special
 cases.

In this paper Gaussian measures will be investigated on the 3--dimensional
 Heisenberg group \ $\HH$ \ which can be obtained by furnishing \ $\RR^3$
 \ with its natural topology and with the product
 $$
  (g_1,g_2,g_3)(h_1,h_2,h_3)
  =\Big(g_1+h_1,g_2+h_2,g_3+h_3+\frac{1}{2}(g_1h_2-g_2h_1)\Big).
 $$
The Schr\"odinger representations \ $\{\pi_{\pm\lambda}:\lambda>0\}$ \ of
 \ $\HH$ \ are representations in the group of unitary operators of the
 complex Hilbert space \ $L^2(\RR)$ \ given by
 \begin{equation}
  \label{SCHREP}
  [\pi_{\pm\lambda}(g)u](x)
  :=\ee^{\pm i(\lambda g_3+\sqrt{\lambda}g_2x+\lambda g_1g_2/2)}
    u(x+\sqrt{\lambda}g_1)
 \end{equation}
 for \ $g=(g_1,g_2,g_3)\in\HH$, \ $u\in L^2(\RR)$ \ and \ $x\in\RR$.
\ The value of the Fourier transform of a probability measure \ $\mu$ \ on
 \ $\HH$ \ at the Schr\"odinger representation \ $\pi_{\pm\lambda}$ \ is the
 bounded linear operator
 \ $\hmu(\pi_{\pm\lambda}):L^2(\RR)\to L^2(\RR)$ \ given by
 $$
  \hmu(\pi_{\pm\lambda})u:=\int_{\HH}\pi_{\pm\lambda}(g)u\,\mu(\dd g),
  \qquad u\in L^2(\RR),
 $$
 interpreted as a Bochner integral.

Let \ $(\mu_t)_{t\sgeq0}$ \ be a Gaussian convolution semigroup of
 probability measures on \ $\HH$ \ (see Section \ref{PREL}).
\ By a result of Siebert \cite[Proposition 3.1, Lemma 3.1]{SIE81},
 \ $\big(\widehat{\mu_t}(\pi_{\pm\lambda})\big)_{t\sgeq0}$ \ is a strongly
 continuous semigroup of contractions on \ $L^2(\RR)$ \ with infinitesimal
 generator
 $$
  N(\pi_{\pm\lambda})
  =\alpha_1I+\alpha_2x+\alpha_3D+\alpha_4x^2+\alpha_5(xD+Dx)+\alpha_6D^2,
 $$
 where \ $\alpha_1,\dots,\alpha_6$ \ are certain complex numbers (depending on
 \ $(\mu_t)_{t\sgeq0}$, \ see Remark \ref{MEHLER}), \ $I$ \ denotes the
 identity operator on \ $L^2(\RR)$, \ $x$ \ is the multiplication by the
 variable \ $x$, \ and \ $Du(x)=u'(x)$.
\ One of our purposes is to determine the action of the operators
 $$
  \widehat{\mu_t}(\pi_{\pm\lambda})=\ee^{tN(\pi_{\pm\lambda})},\qquad t\geq0
 $$
 on \ $L^2(\RR)$. (Here the notation \ $(e^{tA})_{t\sgeq 0}$ \ means
 a semigroup of operators with
 infinitesimal generator \ $A$.)
When \ $N(\pi_{\pm\lambda})$ \ has the special form \ $\frac{1}{2}(D^2-x^2)$,
 \ the celebrated Mehler's formula gives us
 $$
  \ee^{t(D^2-x^2)/2}u(x)
  =\frac{1}{\sqrt{2\pi\sinh t}}
   \int_\RR\exp\left\{-\frac{(x^2+y^2)\cosh t-2xy}{2\sinh t}\right\}u(y)\,\dd y
 $$
 for all \ $t>0$, \ $u\in L^2(\RR)$ \ and \ $x\in\RR$, \ see, e.g., Taylor
 \cite{TAY86}, Davies \cite{DAV89}.
Our Theorem \ref{SREP} in Section \ref{FourierGaus} can be regarded as a
 generalization of Mehler's formula.

It turns out that
 \ $\widehat{\mu_t}(\pi_{\pm\lambda})=\ee^{tN(\pi_{\pm\lambda})}$, \ $t\geq0$
 \ are again integral operators on \ $L^2(\RR)$ \ if \ $\alpha_6$ \ is a
 positive real number.
One of the main results of the present paper is an explicit formula for the
 kernel function of these integral operators (see Theorem \ref{SREP}).
We apply a probabilistic method using that the Fourier transform
 \ $\hmu(\pi_{\pm\lambda})$ \ of an absolutely continuous probability
 measure \ $\mu$ \ on \ $\HH$ \ can be derived from the Euclidean Fourier
 transform of \ $\mu$ \ considering \ $\mu$ \ as a measure on \ $\RR^3$ \ (see
 Proposition \ref{ACP}). We note that a random walk approach might
 provide a different proof of Theorem \ref{SREP}, but we think that
 it would not be simpler than ours.

The second part of the paper deals with convolutions of Gaussian measures on
 \ $\HH$.
\ The convolution of two Gaussian measures on a locally compact Abelian group
 is again a Gaussian measure (it can be proved by the help of Fourier
 transforms; see Parthasarathy \cite{PAR67}).
We prove that a convolution of Gaussian measures on \ $\HH$ \ is almost never a
 Gaussian measure. More exactly, we obtain the following result (using our
 explicit formula for the Fourier transforms).

\begin{Thm}\label{C} \
Let \ $\mu'$ \ and \ $\mu''$ \ be Gaussian measures on \ $\HH$.
\ Then the convolution \ $\mu'*\mu''$ \ is a Gaussian measure on \ $\HH$ \ if
 and only if one of the following conditions holds:
\renewcommand{\labelenumi}{{\rm(C\arabic{enumi})}}
 \begin{enumerate}
  \item there exist elements \ $Y_0',\,Y_0'',\,Y_1,\,Y_2$ \ in the Lie algebra
         of \ $\HH$ \ such that \ $[Y_1,Y_2]=0$, \ the support of \ $\mu'$ \ is
         contained in \ $\exp\{Y_0'+\RR\cdot Y_1+\RR\cdot Y_2\}$ \ and the
         support of \ $\mu''$ \ is contained in
         \ $\exp\{Y_0''+\RR\cdot Y_1+\RR\cdot Y_2\}$.
        \ (Equivalently, there exists an Abelian subgroup \ $\GG$ \ of \ $\HH$
         \ such that \ $\supp(\mu')$ \ and \ $\supp(\mu'')$ \ are contained in
         ``Eucledian cosets'' of \ $\GG$.)
  \item there exist a Gaussian semigroup \ $(\mu_t)_{t\sgeq0}$ \ and
         \ $t',t''\geq0$ \ and a Gaussian measure \ $\nu$ \ such that
         \ $\supp(\nu)$ \ is contained in the center of \ $\HH$ \ and either
         \ $\mu'=\mu_{t'}$, \ $\mu''=\mu_{t''}*\nu$ \ or \ $\mu'=\mu_{t'}*\nu$,
         \ $\mu''=\mu_{t''}$ \ holds.
        \ (Equivalently, \ $\mu'$ \ and \ $\mu''$ \ are sitting on the same
         Gaussian semigroup modulo a Gaussian measure with support contained in
         the center of \ $\HH$. )
 \end{enumerate}
\end{Thm}

We note that in case of \ ${\rm (C1)}$, \ $\mu'$ \ and \ $\mu''$ \ are Gaussian
 measures also in the ``Euclidean sense'' (i.e., considering them as measures
 on \ $\RR^3$).
Moreover, Theorem \ref{CONV} contains an explicit formula for the
Fourier transform of a convolution of arbitrary Gaussian measures on
\ $\HH$. \

The structure of the present work is similar to Pap \cite{PAP02}.
Theorems \ref{C} and \ref{SREP} of the present paper are
generalizations of the corresponding results for symmetric Gaussian
measures on \ $\HH$ \ due to Pap \cite{PAP02}. We summarize briefly
the new ingredients needed in the present paper. Comparing Lemma 6.1
in Pap \cite{PAP02} and Proposition \ref{jfEF} of the present paper,
one can realize that now we have to calculate a much more
complicated (Euclidean) Fourier transform (see \eqref{NINE}). For
this reason we generalized a result due to Chaleyat-Maurel
\cite{CHA81} (see Lemma \ref{CHAMAU}). We note that using Lemma
\ref{GAUSSPAR2} one can easily derive Theorem 1.1 in Pap
\cite{PAP02} from Theorem \ref{C} of the present paper.

It is natural to ask whether we can prove our results for
non-symmetric Gaussian measures using only the results for symmetric
Gaussian measures. The answer is no.  The reason for this is that in
case of \ $\HH$ \ the convolution of a symmetric Gaussian measure
and a Dirac measure is in general not a Gaussian measure. For
example, if \ $a=(1,0,0)\in\HH$ \ and \ $(\mu_t)_{t\sgeq 0}$ \ is a
Gaussian semigroup with infinitesimal generator \ $X_1^2+X_2^2,$ \
then using Lemma \ref{GAUSSSUPPORT}, one can easily check that \
$\mu_1*\varepsilon_a$ \ is not a Gaussian measure on \ $\HH,$ \
where \ $\varepsilon_a$ \ denotes the Dirac measure concentrated on
the element \ $a\in\HH.$ \ (For the definition of an infinitesimal
generator and \ $X_1,X_2,X_3,$ \ see Section 2.)

We note that if the convolution of two Gaussian measures on \ $\HH$
\ is again a Gaussian measure on \ $\HH,$ \ then the corresponding
infinitesimal generators not neccesarily commute, nor even if the
infinitesimal generator corresponding to the convolution is the sum
of the original infinitesimal generators. Now we give an
illuminating counterexample. Let \ $\mu'$ \ and \ $\mu''$ \ be
Gaussian measures on \ $\HH$ \ such that the corresponding Gaussian
semigroups have infinitesimal generators
$$
  N'=\frac{1}{2}(X_1+X_2)^2
   \qquad\text{and}\qquad
  N''=\frac{1}{2}(X_1+X_2)^2+X_1X_3,
  \qquad\text{respectively.}
$$
 Using Theorem \ref{SCHREPPREL} and Lemma \ref{GAUSSPAR2},
 \ $\mu'*\mu''$ \ is a symmetric Gaussian measure on \ $\HH$ \
such that the corresponding Gaussian semigroup has infinitesimal
generator \ $N'+N''.$ \ But \ $N'$ \ and \ $N''$ \ do not commute.
Indeed, \ $N'N''-N''N'=-(X_1+X_2)X_3^2\ne0.$

At the end of our paper we formulate Theorem \ref{C} in the
important special case of centered Gaussian measures for which the
corresponding Gaussian semigroups are stable in the sense of Hazod.
This kind of Gaussian measures arises in a standard version of
central limit theorems on \ $\HH$ \ proved by Wehn \cite{WEH62}. In
this special case Theorem \ref{C} can be derived from the results
for symmetric Gaussian measures in Pap \cite{PAP02}.

\section{Preliminaries}\label{PREL}

The Heisenberg group \ $\HH$ \ is a Lie group with Lie algebra \ $\cH$,
 \ which can be realized as the vector space \ $\RR^3$ \ furnished with
 multiplication
 $$
  [(p_1,p_2,p_3),(q_1,q_2,q_3)]=(0,0,p_1q_2-p_2q_1).
 $$
An element \ $X\in\cH$ \ can be regarded as a left--invariant differential
 operator on \ $\HH$, \ namely, for continuously differentiable functions
 \ $f:\HH\to\RR$ \ we put
 $$
  Xf(g):=\lim_{t\to0}t^{-1}\Big(f(g\exp(tX))-f(g)\Big),\qquad g\in\HH,
 $$
 where the exponential mapping \ $\exp:\cH\to\HH$ \ is now the identity
 mapping.

A family \ $(\mu_t)_{t\sgeq0}$ \ of probability measures on \ $\HH$ \ is said
 to be a \emph{(continuous) convolution semigroup} if we have
 \ $\mu_s*\mu_t=\mu_{s+t}$ \ for all \ $s,t\geq0$, \ and
 \ $\lim_{t\downarrow0}\mu_t=\mu_0=\vare_e$, \ where \ $e=(0,0,0)$ \
 is the unit element of \ $\HH.$ \ Its \emph{infinitesimal generator} is defined by
 $$
  (Nf)(g)
  :=\lim_{t\downarrow0}t^{-1}\int_\HH\big(f(gh)-f(g)\big)\mu_t(\dd h),
  \qquad g\in\HH,
 $$
 for suitable functions \ $f:\HH\to\RR$.
(The infinitesimal generator is always defined for infinitely differentiable
 functions \ $f:\HH\to\RR$ \ with compact support.)
A convolution semigroup \ $(\mu_t)_{t\sgeq0}$ \ is called a
 \emph{Gaussian semigroup} if
 \ $\lim_{t\downarrow0}t^{-1}\mu_t(\HH\setminus U)=0$ \ for all (Borel)
 neighbourhoods \ $U$ \ of \ $e$.
\ Let \ $\{X_1,X_2,X_3\}$ \ denote the natural basis in \ $\cH$ \ (that is,
 \ $\exp X_1=(1,0,0)$, \ $\exp X_2=(0,1,0)$ \ and \ $\exp X_3=(0,0,1)$).
It is known that a convolution semigroup \ $(\mu_t)_{t\sgeq0}$ \ is a
 Gaussian semigroup if and only if its infinitesimal generator has the form
 \begin{equation}
  \label{GIG}
  N=\sum_{k=1}^3a_kX_k+\frac{1}{2}\sum_{j=1}^3\sum_{k=1}^3b_{j,k}X_jX_k,
 \end{equation}
 where \ $a=(a_1,a_2,a_3)\in\RR^3$ \ and \ $B=(b_{j,k})_{1\sleq j,k\sleq3}$
 \ is a real, symmetric, positive semidefinite matrix.
A probability measure \ $\mu$ \ on \ $\HH$ \ is called a
 \emph{Gaussian measure} if there exists a Gaussian semigroup
 \ $(\mu_t)_{t\sgeq0}$ \ such that \ $\mu=\mu_1$.
\ A Gaussian measure on \ $\HH$ \ can be embedded only in a uniquely determined
 Gaussian semigroup (see Baldi \cite{BAL85}, Pap \cite{PAP94}).
(Neuenschwander \cite{NEU93} showed that a Gaussian measure on \ $\HH$ \ can
 not be embedded in a non--Gaussian convolution semigroup.)
Thus for a vector \ $a=(a_1,a_2,a_3)\in\RR^3$ \ and a real, symmetric, positive
 semidefinite matrix \ $B=(b_{j,k})_{1\sleq j,k\sleq3}$ \ we can speak about
 the Gaussian measure \ $\mu$ \ with parameters \ $(a,B)$ \ which is by
 definition \ $\mu:=\mu_1$, \ where \ $(\mu_t)_{t\sgeq0}$ \ is the Gaussian
 semigroup with infinitesimal generator \ $N$ \ given by (\ref{GIG}).
If \ $\nu$ \ is a Gaussian measure with parameters \ $(a,B)$ \ and
 \ $(\nu_s)_{s\sgeq0}$ \ is the Gaussian semigroup with infinitesimal generator
 \ $N$ \ given by (\ref{GIG}) then \ $\nu_t$ \ is a Gaussian measure with
 parameters \ $(ta,tB)$ \ for all \ $t\geq0$, \ since \ $\mu_s:=\nu_{st}$,
 \ $s\geq0$ \ defines a Gaussian semigroup with infinitesimal generator \ $tN$.
\ Hence \ $\nu_t=\mu_1$, \ so it will be sufficient to calculate the Fourier
 transform of \ $\mu_1$.

Let us consider a Gaussian semigroup \ $(\mu_t)_{t\sgeq0}$ \ with parameters
 \ $(a,B)$ \ on \ $\HH$.
\ Its infinitesimal generator \ $N$ \ can be also written in the form
 \begin{equation}
  \label{GIGM}
  N=Y_0+\frac{1}{2}\sum_{j=1}^dY_j^2,
 \end{equation}
 where \ $0\leq d\leq3$ \ and
 $$
  Y_0=\sum_{k=1}^3a_kX_k,\qquad
  Y_j=\sum_{k=1}^3\sigma_{k,j}X_k,\quad 1\leq j\leq d,
 $$
 where \ $\Sigma=(\sigma_{k,j})$ \ is a \ $3\times d$ \ matrix with
 \ $\rank(\Sigma)=\rank(B)=d$.
\ Moreover, \ $B=\Sigma\cdot\Sigma^\top$.
\ Then the measure \ $\mu_t$ \ can be described as the distribution of the
 random vector \ $Z(t)=(Z_1(t),Z_2(t),Z_3(t))$ \ with values in \ $\RR^3$,
 \ where
 \begin{gather*}
  Z_1(t)=a_1t+\sum_{k=1}^d\sigma_{1,k}W_k(t),\qquad
  Z_2(t)=a_2t+\sum_{k=1}^d\sigma_{2,k}W_k(t),\\
  \begin{align*}
   Z_3(t)&=a_3t+\sum_{k=1}^d\sigma_{3,k}W_k(t)
           +\frac{1}{2}\int_0^t(Z_1(s)\,\dd Z_2(s)-Z_2(s)\,\dd Z_1(s))\\
         &=a_3t+\sum_{k=1}^d\sigma_{3,k}W_k(t)
           +\sum_{k=1}^d(a_2\sigma_{1,k}-a_1\sigma_{2,k})W^*_k(t)
           +\sum_{1\sleq k<\ell\sleq d}
             (\sigma_{1,k}\sigma_{2,\ell}-\sigma_{1,\ell}\sigma_{2,k})
             W_{k,\ell}(t),
  \end{align*}
 \end{gather*}
 where \ $(W_1(t),\ldots,W_d(t))_{t\sgeq0}$ \ is a standard Wiener process in
 \ $\RR^d$ \ and
 \begin{align*}
  W^*_k(t)
   &:=\frac{1}{2}\left(\int_0^tW_k(s)\,\dd s-\int_0^ts\,\dd W_k(s)\right),\\
  W_{k,\ell}(t)
   &:=\frac{1}{2}\left(\int_0^tW_k(s)\,\dd W_\ell(s)
                       -\int_0^tW_\ell(s)\,\dd W_k(s)\right).
 \end{align*}
(See, e.g., Roynette \cite{ROY75}.)
The process \ $(W_{k,\ell}(t))_{t\sgeq0}$ \ is the so--called L\'evy's
 stochastic area swept by the process \ $(W_k(s),W_\ell(s))_{s\in[0,t]}$ \ on
 \ $\RR^2$.

\section{Fourier transform of a Gaussian measure}\label{FourierGaus}

The Schr\"odinger representations are infinite dimensional, irreducible,
 unitary representations, and each irreducible, unitary representation is
 unitarily equivalent with one of the Schr\"odinger representations or with
 \ $\chi_{\alpha,\beta}$ \ for some \ $\alpha,\beta\in\RR$, \ where
 \ $\chi_{\alpha,\beta}$ \ is a one--dimensional representation given by
 $$
  \chi_{\alpha,\beta}(g):=\ee^{i(\alpha g_1+\beta g_2)},\qquad
  g=(g_1,g_2,g_3)\in\HH.
 $$
The value of the \emph{Fourier transform} of a probability measure \ $\mu$
 \ on \ $\HH$ \ at the representation \ $\chi_{\alpha,\beta}$ \ is
 $$
  \hmu(\chi_{\alpha,\beta})
  :=\int_{\HH}\chi_{\alpha,\beta}(g)\,\mu(\dd g)
   =\int_{\HH}\ee^{i(\alpha g_1+\beta g_2)}\,\mu(\dd g)
   =\widetilde{\mu}(\alpha,\beta,0),
 $$
 where \ $\widetilde{\mu}:\RR^3\to\CC$ \ denotes the \emph{Euclidean Fourier
 transform} of \ $\mu$,
 $$
  \widetilde{\mu}(\alpha,\beta,\gamma)
  :=\int_{\HH}\ee^{i(\alpha g_1+\beta g_2+\gamma g_3)}\,\mu(\dd g).
 $$
\indent
Let us consider a Gaussian semigroup \ $(\mu_t)_{t\sgeq0}$ \ with parameters
 \ $(a,B)$ \ on \ $\HH$.
\ The Fourier transform of \ $\mu:=\mu_1$ \ at the one--dimensional
 representations can be calculated easily, since the description of
 \ $(\mu_t)_{t\sgeq0}$ \ given in Section \ref{PREL} implies that
 $$
  \hmu(\chi_{\alpha,\beta})
  =\EE\exp\left\{i(\alpha a_1+\beta a_2)
                 +i\left(\alpha\sum_{k=1}^d\sigma_{1,k}W_k(1)
                         +\beta\sum_{k=1}^d\sigma_{2,k}W_k(1)\right)\right\}
 $$
 for \ $\alpha,\beta\in\RR$.
\ The random variable
 $$
  \left(\sum_{k=1}^d\sigma_{1,k}W_k(1),\,\sum_{k=1}^d\sigma_{2,k}W_k(1)\right)
 $$
 has a normal distribution with zero mean and covariance matrix
 $$
  \begin{bmatrix}
   \sigma_{1,1} & \dots & \sigma_{1,d} \\
   \sigma_{2,1} & \dots & \sigma_{2,d}
  \end{bmatrix}
  \begin{bmatrix}
   \sigma_{1,1} & \sigma_{2,1} \\
   \vdots & \vdots \\ \sigma_{1,d} & \sigma_{2,d}
  \end{bmatrix}
   =\begin{bmatrix}
     b_{1,1} & b_{1,2} \\
     b_{2,1} & b_{2,2}
    \end{bmatrix},
 $$
 since \ $\Sigma\Sigma^\top=B$.
\ Consequently,
 $$
  \hmu(\chi_{\alpha,\beta})
   =\exp\left\{i(\alpha a_1+\beta a_2)
               -\frac{1}{2}
               (b_{1,1}\alpha^2+2b_{1,2}\alpha\beta+b_{2,2}\beta^2)\right\}.$$
\indent
One of the main results of the present paper is an explicit formula for the
 Fourier transform of a Gaussian measure on the Heisenberg group \ $\HH$ \ at
 the Schr\"odinger representations.

\begin{Thm}\label{SREP} \
Let \ $\mu$ \ be a Gaussian measure on \ $\HH$ \ with parameters \ $(a,B)$.
\ Then
 $$
  [\hmu(\pi_{\pm\lambda})u](x)
  =\begin{cases}
    \DS\int_{\RR} K_{\pm\lambda}(x,y)u(y)\,\dd y
     & \text{if \ $b_{1,1}>0$,}\\[4mm]
    L_{\pm\lambda}(x)u(x+\sqrt{\lambda}a_1) & \text{if \ $b_{1,1}=0$,}
   \end{cases}
 $$
 for \ $u\in L^2(\RR)$, \ $x\in\RR$, \ where
 $$
  K_{\pm\lambda}(x,y)
  :=C_{\pm\lambda}(B)
    \exp\left\{-\frac{1}{2}
                \mathbf{z}^\top D_{\pm\lambda}(a,B)\mathbf{z}\right\},
  \qquad\mathbf{z}:=(x,y,1)^\top,
 $$
 where, with \ $\delta:=\sqrt{b_{1,1}b_{2,2}-b_{1,2}^2}$,
 \ $\delta_1:=b_{1,1}b_{2,3}-b_{1,2}b_{1,3}$,
 \ $\delta_2:=a_1b_{1,2}-a_2b_{1,1}$,
 $$
  C_{\pm\lambda}(B)
  :=\begin{cases}
     \DS\frac{1}{\sqrt{2\pi\lambda b_{1,1}}}
      & \text{if \ $\delta=0$,} \\[7mm]
     \DS\sqrt{\frac{\delta}{2\pi b_{1,1}\sinh(\lambda\delta)}}
      & \text{if \ $\delta>0$,}
    \end{cases}
 $$
 and \ $D_{\pm\lambda}(a,B)=(d_{j,k}^{\pm\lambda}(a,B))_{1\sleq j,k\sleq3}$
 \ are symmetric matrices defined for \ $b_{1,1}>0$ \ and \ $\delta=0$ \ by
 $$
  d_{1,1}^{\pm\lambda}(a,B):=\frac{\lambda^{-1}\pm ib_{1,2}}{b_{1,1}},\qquad
  d_{1,2}^{\pm\lambda}(a,B):=-\frac{1}{\lambda b_{1,1}},\qquad
  d_{2,2}^{\pm\lambda}(a,B):=\frac{\lambda^{-1}\mp ib_{1,2}}{b_{1,1}},
 $$
 $$
  d_{1,3}^{\pm\lambda}(a,B)
  :=\frac{a_1\pm i\lambda b_{1,3}}{\sqrt{\lambda}b_{1,1}}
    \pm i\frac{\sqrt{\lambda}\delta_2}{2b_{1,1}},\qquad
  d_{2,3}^{\pm\lambda}(a,B)
  :=-\frac{a_1\pm i\lambda b_{1,3}}{\sqrt{\lambda}b_{1,1}}
     \pm i\frac{\sqrt{\lambda}\delta_2}{2b_{1,1}},
 $$
 $$
  d_{3,3}^{\pm\lambda}(a,B)
  :=\frac{(a_1\pm i\lambda b_{1,3})^2}{b_{1,1}}
    +\frac{\lambda^2\delta_2^2}{12b_{1,1}}
    +\lambda^2b_{3,3}\mp2i\lambda a_3,
 $$
 and for \ $\delta>0$ \ by
 $$
  d_{1,1}^{\pm\lambda}(a,B)
  :=\frac{\delta\coth(\lambda\delta)\pm ib_{1,2}}{b_{1,1}},\quad
  d_{1,2}^{\pm\lambda}(a,B):=-\frac{\delta}{b_{1,1}\sinh(\lambda\delta)},\quad
  d_{2,2}^{\pm\lambda}(a,B)
  :=\frac{\delta\coth(\lambda\delta)\mp ib_{1,2}}{b_{1,1}},
 $$
 $$
  d_{1,3}^{\pm\lambda}(a,B)
  :=\frac{a_1\pm i\lambda b_{1,3}}{\sqrt{\lambda}b_{1,1}}
    +\frac{\lambda\delta_1\pm i\delta_2}
          {\sqrt{\lambda}b_{1,1}\delta\coth(\lambda\delta/2)},\quad
  d_{2,3}^{\pm\lambda}(a,B)
  :=-\frac{a_1\pm i\lambda b_{1,3}}{\sqrt{\lambda}b_{1,1}}
    +\frac{\lambda\delta_1\pm i\delta_2}
          {\sqrt{\lambda}b_{1,1}\delta\coth(\lambda\delta/2)},
 $$
 $$
  d_{3,3}^{\pm\lambda}(a,B)
  :=\frac{(a_1\pm i\lambda b_{1,3})^2}{b_{1,1}}
    -\frac{(\lambda\delta_1\pm i\delta_2)^2}{\lambda b_{1,1}\delta^3}
     \Big(\lambda\delta-2\tanh(\lambda\delta/2)\Big)
    +\lambda^2b_{3,3}\mp2i\lambda a_3,
 $$
 and
 \begin{align*}
  L_{\pm\lambda}(x)
  :=\exp\bigg\{&\pm\frac{i\sqrt{\lambda}}{2}
                \Big(\sqrt{\lambda}(2a_3+a_1a_2)+2a_2x\Big)
                -\frac{\lambda^2}{6}(3b_{3,3}+3a_1b_{2,3}+a_1^2b_{2,2})\\[2mm]
               &-\frac{\lambda^{3/2}}{2}(2b_{2,3}+a_1b_{2,2})x
                -\frac{\lambda}{2}b_{2,2}x^2\bigg\}.
 \end{align*}
\end{Thm}

We prove this theorem in Section \ref{Provfourier}.

\begin{Rem}\label{MEHLER} \
Consider a Gaussian convolution semigroup \ $(\mu_t)_{t\sgeq0}$ \ with
 infinitesimal generator \ $N$ \ given in \eqref{GIG}.
Siebert \cite[Proposition 3.1, Lemma 3.1]{SIE81} proved that
 \ $\big(\widehat{\mu_t}(\pi_{\pm\lambda})\big)_{t\sgeq0}$ \ is a strongly
 continuous semigroup of contractions on \ $L^2(\RR)$ \ with infinitesimal
 generator
 $$
  N(\pi_{\pm\lambda})
  =\sum_{k=1}^3a_kX_k(\pi_{\pm\lambda})
    +\frac{1}{2}
     \sum_{j=1}^3\sum_{k=1}^3b_{j,k}X_j(\pi_{\pm\lambda})X_k(\pi_{\pm\lambda}),
 $$
 where
 $$
  X(\pi_{\pm\lambda})u
  :=\lim_{t\to0}t^{-1}\big(\pi_{\pm\lambda}(\exp(tX))u-u\big)
 $$
for all differentiable vectors \ $u$.
\ Hence
 \begin{align*}
  [X_1(\pi_{\pm\lambda})u](x)&=\sqrt{\lambda}u'(x)=\sqrt{\lambda}Du(x),\\[1mm]
  [X_2(\pi_{\pm\lambda})u](x)&=\pm i\sqrt{\lambda}xu(x),\\[1mm]
  [X_3(\pi_{\pm\lambda})u](x)&=\pm i\lambda u(x)
 \end{align*}
 for all \ $x\in\RR$.
\ Consequently,
 $$
  N(\pi_{\pm\lambda})
  =\alpha_1I+\alpha_2x+\alpha_3D+\alpha_4x^2+\alpha_5(xD+Dx)+\alpha_6D^2,
 $$
 where
 $$
  \alpha_1=-\frac{1}{2}\lambda^2b_{3,3}\pm i\lambda a_3,\qquad
  \alpha_2=-\lambda^{3/2}b_{2,3}\pm i\lambda^{1/2}a_2,\qquad
  \alpha_3=\lambda^{1/2}a_1\pm i\lambda^{3/2}b_{1,3},
 $$
 $$
  \alpha_4=-\frac{1}{2}\lambda b_{2,2},\qquad
  \alpha_5=\pm\frac{i}{2}\lambda b_{1,2},\qquad
  \alpha_6=\frac{1}{2}\lambda b_{1,1}.
 $$
\end{Rem}

\section{Absolute continuity and singularity of Gaussian measures}

A probability measure \ $\mu$ \ on \ $\HH$ \ is said to be
absolutely continuous or singular if it is absolutely continuous or
singular with respect to a (and then necessarily to any) Haar
measure on \ $\HH.$ \ It is known that the class of left Haar
measures on \ $\HH$ \ is the same as the class of right Haar
measures on \ $\HH$ \ and hence we can use the expression ''a Haar
measure on \ $\HH$''. \ It is also known that a measure \ $\nu$ \
 on \ $\HH$ \ is a Haar measure if and only if \ $\nu$ \ is the
 Lebesgue measure on \ $\RR^3$ \ multiplied by some positive
 constant. The following proposition is
 the same as Proposition 2.1 in Pap \cite{PAP02}. But the proof
 given here is simpler, we do not use Weyl calculus.

\begin{Pro}\label{ACP} \
If \ $\mu$ \ is absolutely continuous with density \ $f$ \ then the Fourier
 transform \ $\hmu(\pi_{\pm\lambda})$ \ is an integral operator on
 \ $L^2(\RR)$,
 $$
  [\hmu(\pi_{\pm\lambda})u](x)=\int_\RR K_{\pm\lambda}(x,y)u(y)\,\dd y
 $$
 with kernel function \ $K_{\pm\lambda}:\RR^2\to\CC$ \ given by
 $$
  K_{\pm\lambda}(x,y)
  :=\frac{1}{\sqrt{\lambda}}
    \widetilde{f}_{2,3}
     \left(\frac{y-x}{\sqrt{\lambda}},
           \pm\sqrt{\lambda}\left(\frac{y+x}{2}\right),
           \pm\lambda\right),
 $$
 where
 $$
  \widetilde{f}_{2,3}(s_1,\ts_2,\ts_3)
  :=\int_{\RR^2}\ee^{i(\ts_2s_2+\ts_3s_3)}f(s_1,s_2,s_3)\,\dd s_2\,\dd s_3,
  \qquad(s_1,\ts_2,\ts_3)\in\RR^3
 $$
 denotes a partial Euclidean Fourier transform of \ $f$.
\end{Pro}

\noindent
\textbf{Proof.} \
Using the definition of the Schr\"odinger representation we obtain
 \begin{align*}
  [\hmu(\pi_{\pm\lambda})u](x)
   &=\int_{\RR^3}
      \ee^{\pm i(\lambda s_3+\sqrt{\lambda}s_2x+\lambda s_1s_2/2)}
      u(x+\sqrt{\lambda}s_1)f(s_1,s_2,s_3)\,\dd s_1\,\dd s_2\,\dd s_3\\
   &=\frac{1}{\sqrt{\lambda}}
     \int_{\RR^3}
      \ee^{\pm i(\lambda s_3+\sqrt{\lambda}s_2x+\sqrt{\lambda}(y-x)s_2/2)}
      u(y)f\left(\frac{y-x}{\sqrt{\lambda}},s_2,s_3\right)
      \,\dd y\,\dd s_2\,\dd s_3\\
   &=\int_\RR K_{\pm\lambda}(x,y)u(y)\,\dd y,
 \end{align*}
 where
 \begin{align*}
  K_{\pm\lambda}(x,y)
   &=\frac{1}{\sqrt{\lambda}}
     \int_{\RR^2}
      \ee^{\pm i(\lambda s_3+\sqrt{\lambda}(x+y)s_2/2)}
      f\left(\frac{y-x}{\sqrt{\lambda}},s_2,s_3\right)\,\dd s_2\,\dd s_3\\
   &=\frac{1}{\sqrt{\lambda}}
      \widetilde{f}_{2,3}
       \left(\frac{y-x}{\sqrt{\lambda}},
             \pm\sqrt{\lambda}\left(\frac{y+x}{2}\right),
             \pm\lambda\right).
 \end{align*}
Hence the assertion.
\proofend

The partial Euclidean Fourier transform \ $\tf_{2,3}$ \ can be obtained by the
 inverse Euclidean Fourier transform:
 \begin{equation}
  \label{partialfourier}
  \tf_{2,3}(s_1,\ts_2,\ts_3)
    =\frac{1}{2\pi}\int_\RR\ee^{-is_1\ts_1}
      \tf(\ts_1,\ts_2,\ts_3)\,\dd\ts_1,\qquad
   (s_1,\ts_2,\ts_3)\in\RR^3,
 \end{equation}
 where \ $\tf$ \ denotes the (full) Euclidean Fourier transform of \ $f$\,:
 $$
  \tf(\ts_1,\ts_2,\ts_3)
  :=\int_{\RR^3}\ee^{i(\ts_1 s_1+\ts_2 s_2+\ts_3 s_3)}f(s_1,s_2,s_3)\,
     \dd s_1\,\dd s_2\,\dd s_3
 $$
 for \ $(\ts_1,\ts_2,\ts_3)\in\RR^3$.
\ Moreover, \ $\hmu(\pi_{\pm\lambda})$ \ is a compact operator.
If the density \ $f$ \ of \ $\mu$ \ belongs to the Schwarz space then
 \ $\hmu(\pi_{\pm\lambda})$ \ is a trace class (i.e., nuclear) operator.

In order to apply Proposition \ref{ACP} we shall need the description of the
 set of absolutely continuous Gaussian measures on \ $\HH$.
Using a general result due to Siebert \cite[Theorem 2]{SIE82} one can prove
 the following lemma as in Pap \cite[Lemma 3.3]{PAP02}.

\begin{Lem} \
A Gaussian measure \ $\mu$ \ on \ $\HH$ \ with parameters \ $(a,B)$ \ is
 either absolutely continuous or singular.
More precisely, \ $\mu$ \ is absolutely continuous if
 \ $b_{1,1}b_{2,2}-b_{1,2}^2>0$ \ and singular if
 \ $b_{1,1}b_{2,2}-b_{1,2}^2=0$.
\end{Lem}

The next lemma describes the support of a Gaussian measure on \ $\HH$.

\begin{Lem}\label{GAUSSSUPPORT} \
Let \ $(\mu_t)_{t\sgeq0}$ \ be a Gaussian semigroup on \ $\HH$ \ with
 infinitesimal generator \ $N$ \ given by \eqref{GIGM}.
According to the structure of \ $N$ \ we can distinguish five different types
 of Gaussian semigroups:
 \renewcommand{\labelenumi}{{\rm(\roman{enumi})}}
 \begin{enumerate}
  \item $N=Y_0+\frac{1}{2}(Y_1^2+Y_2^2+Y_3^2)$ \ with \ $Y_1$, \ $Y_2$ \ and
         \ $Y_3$ \ linearly independent.
        Then the semigroup is absolutely continuous and \ $\supp(\mu_t)=\HH$
         \ for all \ $t>0$.
        \ Moreover, \ $\rank(B)=3$, \ $b_{1,1}b_{2,2}-b_{1,2}^2\ne 0$.
  \item $N=Y_0+\frac{1}{2}(Y_1^2+Y_2^2)$ \ with \ $Y_1$ \ and \ $Y_2$
         \ linearly independent and \ $[Y_1,Y_2]\not=0$.
        \ Then the semigroup is absolutely continuous and \ $\supp(\mu_t)=\HH$
         \ for all \ $t>0$.
        Moreover, \ $\rank(B)=2$, \ $b_{1,1}b_{2,2}-b_{1,2}^2\ne 0$.
  \item $N=Y_0+\frac{1}{2}(Y_1^2+Y_2^2)$ \ with \ $Y_1$ \ and \ $Y_2$
         \ linearly independent and \ $[Y_1,Y_2]=0$.
        \ Then the semigroup is singular, it is a Gaussian semigroup on
         \ $\RR^3$ \ as well, and it is supported by a `Euclidean coset' of the
         same closed normal subgroup, namely,
         $$
          \supp(\mu_t)=\exp(tY_0+\RR\cdot Y_1+\RR\cdot Y_2)
         $$
         for all \ $t>0$.
        \ Moreover, \ $\rank(B)=2$, \ $b_{1,1}b_{2,2}-b_{1,2}^2=0$.
  \item $N=Y_0+\frac{1}{2}Y_1^2$.
        \ Then the semigroup is singular, it is a Gaussian semigroup on
         \ $\RR^3$ \ as well, and it is supported by a ``Euclidean coset'' of
         the same closed normal subgroup, namely,
         $$
          \supp(\mu_t)=\exp(tY_0+\RR\cdot Y_1+\RR\cdot [Y_0,Y_1])
         $$
         for all \ $t>0$.
        \ Moreover, \ $\rank(B)=1$, \ $b_{1,1}b_{2,2}-b_{1,2}^2=0$.
  \item $N=Y_0$.
        \ Then the semigroup is singular and consists of point measures,
         namely, \ $\mu_t=\vare_{\exp(tY_0)}$ \ for all \ $t\geq 0$.
\end{enumerate}
\end{Lem}

\noindent
\textbf{Proof.} \
From general results due to Siebert \cite[Theorem 2 and Theorem 4]{SIE82}, it
 follows that a Gaussian measure \ $\mu$ \ on \ $\HH$ \ is absolutely
 continuous if and only if
 \ $\cG:=\cL(Y_i,[Y_j,Y_k]:1\leq i\leq d,\,0\leq j\leq k\leq d)=\RR^3$, \ where
 \ $\cL(\cdot)$ \ denotes the linear hull of the given vectors, and
 \ $Y_i\in\cH$, $0\leq i\leq d$ \ are described in (\ref{GIGM}).
Moreover, the support of \ $\mu_t$ \ is
 $$
  \supp(\mu_t)
   =\overline{\bigcup_{n=1}^\infty\bigg(M\exp\Big(\frac{tY_0}{n}\Big)\bigg)^n}
   \qquad\hbox{for\, all}\qquad t>0,
 $$
 where \ $M$ \ is the analytic subgroup of \ $\HH$ \ corresponding to the Lie
 subalgebra generated by \ $\{Y_i:1\leq i\leq r\}$ \ and the bar denotes the
 closure in \ $\HH$.
\ Clearly \ $[Y_i,Y_i]=0$,
 \ $[Y_i,Y_j]=(\sigma_{1,i}\sigma_{2,j}-\sigma_{1,j}\sigma_{2,i})\,X_3$ \ for
 \ $1\leq i<j\leq d$ \ and \ $[Y,Z]\in\cL(X_3)$ \ for all \ $Y,Z\in\cH$.

We prove only the cases \ ${\rm (iii)}$ \ and \ ${\rm (iv)}$, \ the other cases
 can be proved similarly.

In case of \ ${\rm(iii)}$ \ we have \
$\cG=\cL(Y_1,Y_2,[Y_0,Y_1],[Y_0,Y_2])$. \
  Since \ $[Y_1,Y_2]=0$, \ we have
 \ $\sigma_{1,1}\sigma_{2,2}-\sigma_{1,2}^2=0$, \ so \ $Y_1$ \ and \ $Y_2$
 \ are linearly dependent in their first two coordinates, thus their linear
 independence yields \ $X_3\in\cL(Y_1,Y_2)$. \ Moreover,
 \ $[Y_0,Y_1],[Y_0,Y_2]\in\cL(X_3)\subset\cL(Y_1,Y_2).$ \
 So \ $\cG=\cL(Y_1,Y_2)\not=\RR^3$, \ i.e., the semigroup
 \ $(\mu_t)_{t\sgeq0}$ \ is singular.

To obtain the formula for the support of \ $\mu_t$ \ it is
sufficient to prove
 that
 \ $\Big(M\exp\big(\frac{t}{n}Y_0\big)\Big)^n
    =\exp(tY_0+\RR\cdot Y_1+\RR\cdot Y_2)$
 \ for all \ $t>0$ \ and \ $n\in\NN$, \ where now
 \ $M=\exp(\RR\cdot Y_1+\RR\cdot Y_2)$.
\ The multiplication in \ $\HH$ \ can be reconstructed by the help of the
 Campbell--Haussdorf formula
 $$
  \exp(X)\exp(Y)=\exp\Big(X+Y+\frac{1}{2}[X,Y]\Big),\quad X,Y\in\cH.
 $$
Applying induction by \ $n$ \ gives the assertion.
Indeed, for \ $n=1$ \ we have
 \ $M\exp\left(tY_0\right)=\exp(\RR\cdot Y_1+\RR\cdot Y_2)\exp(tY_0)
    =\exp(tY_0+\RR\cdot Y_1+\RR\cdot Y_2)$,
 \ since \ $X_3\in\cL(Y_1,Y_2)$.
\ Suppose that
 \ $\Big(M\exp\big(\frac{t}{n-1}Y_0\big)\Big)^{n-1}
    =\exp(tY_0+\RR\cdot Y_1+\RR\cdot Y_2)$
 \ holds.
Using the Campbell--Haussdorf formula and the induction hypothesis we get
 \ $\Big(M\exp\big(\frac{t}{n}Y_0\big)\Big)^{n}
    =\exp\big(\frac{n-1}{n}tY_0+\RR\cdot Y_1+\RR\cdot Y_2\big)
     \exp\big(\frac{t}{n}Y_0+\RR\cdot Y_1+\RR\cdot Y_2\big)$.
\ Since \ $X_3\in\cL(Y_1,Y_2)$ \ and \ $[Y,Z]\in\cL(X_3)$ \ for all
 \ $Y,\,Z\in\cH$, \ application of the Campbell--Haussdorf formula once more
 gives the assertion.

The case \ ${\rm (iv)}$ \ can be obtained similarly.
Indeed, we have \ $\cG=\cL(Y_1,[Y_0,Y_1])\ne\RR^3$, \ $M=\exp(\RR\cdot Y_1)$,
 \ hence \ $\supp(\mu_t)=\exp(tY_0+\RR\cdot Y_1+\RR\cdot[Y_1,Y_0])$ \ for all
 \ $t>0$.
\proofend

\section{Euclidean Fourier transform of a Gaussian measure and the proof of
Theorem \ref{SREP}}\label{Provfourier}

Now we investigate the processes \ $(W^*_k(t))_{t\sgeq0}$ \ and
 \ $(W_{k,\ell}(t))_{t \sgeq0}$ \ (defined in Section \ref{PREL}).
Let \ $t>0$ \ be fixed.
We prove that \ $W^*_k(t)$ \ and \ $W_{k,\ell}(t)$ \ can be constructed by the
 help of infinitely many independent identically distributed real random
 variables with standard normal distribution.
Because of the self--similarity property of the Wiener process it is sufficient
 to prove the case \ $t=2\pi$.

\begin{Lem}\label{CHAMAUPREL} \
Let \ $(W_1(s),\ldots,W_d(s))_{s\in[0,2\pi]}$ \ be a standard Wiener process
 in \ $\RR^d$ \ on a probability space \ $(\Omega,\cA,P)$.
\ Let us consider the orthonormal basis
 \ $f_n(s)=(2\pi)^{-1/2}e^{ins},\quad s\in[0,2\pi],\quad n\in\ZZ$ \ in the
 complex Hilbert space \ $L^2([0,2\pi])$.
\ If \ $(g(s))_{s\in[0,2\pi]}$ \ is an adapted, measurable, complex valued
 process, independent of \ $(W_1(s),\ldots,W_d(s))_{s\in [0,2\pi]}$ \ such
 that \ $\EE\left(\int_0^{2\pi}\vert g(s)\vert^2\,\dd s\right)<\infty$ \ then
 \begin{equation}
  \label{WCON3}
  \int_0^{2\pi}g(s)\,\dd W_j(s)
  =\sum_{n\in\ZZ}\langle g,f_n\rangle\int_{0}^{2\pi}f_n(s)\,\dd W_j(s)
                               \quad\mathrm{a.\,s.},\qquad j=1,\ldots,d,
 \end{equation}
 where \ $\langle \cdot\,,\cdot\rangle$ \ denotes the inner product in
 \ $L^2([0,2\pi])$ \ and the convergence of the series on the right hand side
 of \eqref{WCON3} is meant in \ $L^2(\Omega,\cA,P).$
\end{Lem}

\noindent \textbf{Proof.} \
Let \ $1\leq j\leq d$ \ be arbitrary, but fixed.
First we prove that the right hand side of \eqref{WCON3} is
 convergent in \ $L^2(\Omega,\cA,P).$ \ Using that the processes \ $(g(s))_{s\in[0,2\pi]}$ \
 and \ $(W_1(s),\ldots,W_d(s))_{s\in [0,2\pi]}$ \ are independent
 of each other, for \ $n,m\in\ZZ,$ \ $n\ne m,$ \ we get
\begin{align*}
\EE\bigg(\langle g,f_n\rangle\int_{0}^{2\pi}f_n(s)\,\dd W_j(s)\,
 &\overline{\langle g,f_m\rangle}\int_{0}^{2\pi}\overline{f_m(s)}\,\dd W_j(s)\bigg)\\
 &=\EE\big(\langle g,f_n\rangle\overline{\langle g,f_m\rangle}\big)
     \EE\bigg(\int_{0}^{2\pi}f_n(s)\,\dd W_j(s)
        \int_{0}^{2\pi}\overline{f_m(s)}\,\dd W_j(s)\bigg)\\
 &=\EE\big(\langle g,f_n\rangle\overline{\langle g,f_m\rangle}\big)
         \int_{0}^{2\pi}f_n(s)\overline{f_m(s)}\,\dd s=0.
\end{align*}
Using again the independence of \ $(g(s))_{s\in[0,2\pi]}$ \
 and \ $(W_1(s),\ldots,W_d(s))_{s\in [0,2\pi]},$ \ we have
\begin{align*}
 \EE\bigg\vert\langle g,f_n\rangle\int_{0}^{2\pi}f_n(s)
                                              \,\dd W_j(s)\bigg\vert^2
  &=\EE\big\vert\langle g,f_n\rangle\big\vert^2
       \EE\bigg\vert\int_{0}^{2\pi}f_n(s)\,\dd W_j(s)\bigg\vert^2\\
  &=\EE\big\vert\langle g,f_n\rangle\big\vert^2
               \int_{0}^{2\pi}\vert f_n(s)\vert^2\,\dd s
  =\EE\big\vert\langle g,f_n\rangle\big\vert^2.
 \end{align*}
Since \ $\EE\left(\int_0^{2\pi}\vert g(s)\vert^2\,\dd
 s\right)<\infty,$ \
 Parseval's identity in \ $L^2([0,2\pi])$ \ gives us that
\begin{align*}
\sum_{n\in\ZZ}\big\vert\langle g,f_n\rangle\big\vert^2
  =\int_{0}^{2\pi}\vert g(s)\vert^2\,\dd s\quad\mathrm{a.\,s.}
\end{align*}
This implies that
\begin{align*}
\sum_{n\in\ZZ}\EE\big\vert\langle g,f_n\rangle\big\vert^2
  =\EE\int_{0}^{2\pi}\vert g(s)\vert^2\,\dd s<\infty.
\end{align*}
Hence the right hand side of \eqref{WCON3} is
 convergent in \ $L^2(\Omega,\cA,P).$ \

We show now that
\begin{align*}
\EE\bigg\vert
     \int_0^{2\pi}g(s)\,\dd W_j(s)
       -\sum_{n\in\ZZ}\langle g,f_n\rangle\int_{0}^{2\pi}f_n(s)\,\dd W_j(s)
   \bigg\vert^2=0,
\end{align*}
which implies \eqref{WCON3}. \ We have
\begin{align*}
 &\EE\bigg\vert
     \int_0^{2\pi}g(s)\,\dd W_j(s)
       -\sum_{n\in\ZZ}\langle g,f_n\rangle\int_{0}^{2\pi}f_n(s)\,\dd W_j(s)
   \bigg\vert^2\\
 &=\EE\bigg\vert\int_{0}^{2\pi}g(s)\,\dd W_j(s)\bigg\vert^2
   +\EE\bigg\vert\sum_{n\in\ZZ}\langle g,f_n\rangle
                   \int_{0}^{2\pi}f_n(s)\,\dd W_j(s)\bigg\vert^2\\
 &\phantom{=\:}
  -2\,\Re\EE\bigg(\int_{0}^{2\pi}g(s)\,\dd W_j(s)
           \sum_{n\in\ZZ}\overline{\langle g,f_n\rangle}
                   \int_{0}^{2\pi}\overline{f_n(s)}\,\dd W_j(s)\bigg)
  =:A_1+A_2-2\,\Re A_3.
\end{align*}
Then we get
\begin{align*}
 A_1&=\EE\int_{0}^{2\pi}\vert g(s)\vert^2\,\dd s,\\
 A_2&=\sum_{n\in\ZZ}\EE\bigg\vert\langle g,f_n\rangle
                           \int_{0}^{2\pi}f_n(s)\,\dd W_j(s)\bigg\vert^2
     =\sum_{n\in\ZZ}\EE\big\vert\langle g,f_n\rangle\big\vert^2
     =\EE\int_{0}^{2\pi}\vert g(s)\vert^2\,\dd s,\\
 A_3&=\sum_{n\in\ZZ}\EE\bigg(\int_{0}^{2\pi}g(s)\,\dd W_j(s)\,
    \overline{\langle g,f_n\rangle}\int_{0}^{2\pi}\overline{f_n(s)}\,\dd W_j(s)\bigg).
\end{align*}
Let us denote the \ $\sigma$-algebra generated by the process \
$(g(s))_{s\in[0,2\pi]}$ \ by \ $\cF(g).$ \ Then we obtain
\begin{align*}
A_3&=\sum_{n\in\ZZ}\EE\,\EE\bigg(\int_{0}^{2\pi}g(s)\,\dd W_j(s)\,
    \overline{\langle g,f_n\rangle}\int_{0}^{2\pi}\overline{f_n(s)}\,\dd W_j(s)\Big\vert
        \,\cF(g)\bigg)\\
    &=\sum_{n\in\ZZ}\EE\bigg(\overline{\langle g,f_n\rangle}\,
             \EE\bigg(\int_{0}^{2\pi}g(s)\,\dd W_j(s)
                        \int_{0}^{2\pi}\overline{f_n(s)}\,\dd
                        W_j(s)\Big\vert\,\cF(g)\bigg)\bigg)\\
    &=\sum_{n\in\ZZ}\EE\bigg(\overline{\langle g,f_n\rangle}
               \int_{0}^{2\pi}g(s)\overline{f_n(s)}\,\dd s\bigg)
    =\sum_{n\in\ZZ}\EE\big\vert\langle g,f_n\rangle\big\vert^2
    =\EE\int_{0}^{2\pi}\vert g(s)\vert^2\,\dd s.
\end{align*}
Hence the assertion.
 \proofend

The next statement is a generalization of Section 1.2 in
Chaleyat-Maurel \cite{CHA81}.

\begin{Lem}\label{CHAMAU} \
Let \ $(W_1(s),\ldots,W_d(s))_{s\in[0,2\pi]}$ \ be a standard Wiener process in
 \ $\RR^d$.
\ Then there exist random variables \ $a_n^{(j)},\,b_n^{(j)}$,
 \ $n\in\NN,\,j=1,\ldots,d$ \ with standard normal distribution, independent of
 each other and of the random variable \ $(W_1(2\pi),\ldots,W_d(2\pi))$ \ such
 that the following constructions hold
 \begin{align}\label{WCON1}
  W_{j,k}(2\pi)
  &=\sum_{n=1}^\infty\frac{1}{n}
     \bigg[b_n^{(j)}\bigg(a_n^{(k)}-\frac{1}{\sqrt{\pi}}W_k(2\pi)\bigg)
           -b_n^{(k)}\bigg(a_n^{(j)}-\frac{1}{\sqrt{\pi}}W_j(2\pi)\bigg)\bigg]
  \qquad{\rm a.\,s.,}\\\label{WCON2}
  W_\ell^*(2\pi)&=-2\sqrt{\pi}\sum_{n=1}^\infty\frac{b_n^{(\ell)}}{n}
  \qquad{\rm a.\,s.}
 \end{align}
 for all \ $1\leq j<k\leq d$ \ and \ $\ell=1,\ldots,d$, \ where the
 series on the right hand sides of \ $(5.2)$ \ and \ $(5.3)$ \ are
 convergent almost surely.
\end{Lem}

\noindent \textbf{Proof.} \
 Retain the notations of Lemma \ref{CHAMAUPREL} and let us denote
 \ $c_n^{(j)}:=\int_0^{2\pi}\overline{f_n(s)}\,\dd W_j(s)$, \ $n\in\ZZ$,
 \ $j=1,\ldots,d.$ \ Then \ $c_n^{(j)}$, \ $n\in\ZZ$, $n\ne 0,$
 \ $j=1,\ldots,d$ \ are independent identically distributed complex random
 variables with standard normal distribution, \ i.e., the decompositions
 \ $c_n^{(j)}=\frac{a_n^{(j)}+i\,b_n^{(j)}}{\sqrt{2}}$, \ $n\in\ZZ$,
 $n\ne 0,$ $j=1,\ldots,d$ \ hold with independent identically distributed real random
 variables \ $a_n^{(j)},\,b_n^{(j)}$, \ $n\in\ZZ$, $n\ne 0,$ $j=1,\ldots,d,$ \ having
 standard normal distribution.
Specifying \ $g$ \ as the indicator function \ $\bone_{[0,t]}$ \ of the
 interval \ $[0,t]$ \ $(t\in[0,2\pi])$ \ in Lemma \ref{CHAMAUPREL}, \
 we have for all \ $t\in[0,2\pi]$
 \begin{equation}
  \label{WCON4}
  W_{\ell}(t)
  =\sum_{n\in\ZZ,\,n\not=0}c_{-n}^{(\ell)}\frac{i}{n}\big(f_{-n}(t)-f_0(t)\big)
   +\frac{c_0^{(\ell)}t}{\sqrt{2\pi}}\quad{\rm a.\,s.,}\qquad \ell=1,\ldots,d.
 \end{equation}
In fact, there is a set \ $\Omega_0$ \ with \ $P(\Omega_0)=0$ \ such that
 \eqref{WCON4} holds for all \ $\omega\notin\Omega_0$ \ and for almost every
 \ $t\in[0,2\pi]$ \ (see, e.g., Ash \cite[p. 107, Problem 4]{ASH72}).
Applying \eqref{WCON3} with \ $g=W_k$ \ and the construction \eqref{WCON4},
 Chaleyat--Maurel \cite{CHA81} showed that \eqref{WCON1} holds.
Choosing \ $g(s)=s\bone_{[0,t]}(s)$, \ $(t\in[0,2\pi])$ \ in Lemma
\ref{CHAMAUPREL} it can be easily checked that
 $$
  \int_0^ts\,\dd W_\ell(s)
  =\sum_{n\in\ZZ,\,n\not=0}\frac{c_{-n}^{(\ell)}(i\,nt+1)}{n^2}f_{-n}(t)
   -\sum_{n\in\ZZ,\,n\not=0}\frac{c_{-n}^{(\ell)}}{n^2}\,\,f_0(t)
   +c_0^{(\ell)}\frac{t^2}{2\sqrt{2\pi}}\qquad{\rm a.\,s.}
 $$
By It\^o's formula we get
 \ $ W_\ell^*(t)=\frac{1}{2}tW_{\ell}(t)-\int_0^ts\,\dd W_{\ell}(s)$.
\ Using the construction \eqref{WCON4} of \ $W_{\ell}(t)$ \ and the definition
 of \ $c_n^{(\ell)}$ \ a simple computation shows that \eqref{WCON2} holds.
By Lemma \ref{CHAMAUPREL} the series in the constructions
\eqref{WCON1}, \eqref{WCON2} and \eqref{WCON4} are convergent in \
$L^2(\Omega,\cA,P)$. \ Since the summands in the series in
\eqref{WCON2} and \eqref{WCON4} are independent, L\'evy's theorem
implies that they are convergent almost surely as well. Finally we
show that the series in \eqref{WCON1} is also convergent almost
surely. For this, using that \ $\sum_{n=1}^\infty b_n^{(\ell)}/n$ \
is convergent almost surely for all \ $\ell=1,\ldots,d,$ \ it is
enough to prove that the series
\begin{align}\label{WCON5}
  \sum_{n=1}^\infty\frac{1}{n}\big(b_n^{(j)}a_n^{(k)}-b_n^{(k)}a_n^{(j)}\big)
\end{align}
 is convergent almost surely. Here \ $b_n^{(j)}a_n^{(k)}-b_n^{(k)}a_n^{(j)},$
 \ $n\in\NN,$ \ are independent, identically distributed real valued
 random variables with zero mean and finite second moment. Hence
 Kolmogorov's One-Series Theorem yields that the series in \eqref{WCON5} is
 convergent almost surely. \proofend

Taking into account Proposition \ref{ACP} and the representation of a Gaussian
 semigroup \ $(\mu_t)_{t\sgeq0}$ \ by the process \ $(Z(t))_{t\sgeq0}$
 \ (given in Section \ref{PREL}), in order to prove Theorem \ref{SREP} we need
 the joint (Euclidean) Fourier transform of the 9--dimensional random vector
 \begin{equation}\label{NINE}
  \big(W_1(t),W_2(t),W_3(t),W_1^*(t),W_2^*(t),W_3^*(t),
       W_{1,2}(t),W_{1,3}(t),W_{2,3}(t)\big).
 \end{equation}

\begin{Pro}\label{jfEF} \
The Fourier transform \ $\tF_t:\RR^9\to\CC$ \ of the random vector \eqref{NINE}
 is
 \begin{align*}
  &\tF_t(\eta_1,\eta_2,\eta_3,\zeta_1,\zeta_2,\zeta_3,
         \xi_{1,2},\xi_{1,3},\xi_{2,3})
   =\\
  &\frac{1}{\cosh(t\|\txi\|/2)}
   \exp\left\{-\frac{t^3}{4\|\txi\|^2}
               \left(\frac{1}{6}-\frac{2\kappa}{t^2\|\txi\|^2}\right)
               \langle\txi,\zeta\rangle^2
              +\frac{\|\txi\|^2\|\teta\|^2
                     +\kappa\langle\txi,\teta\rangle^2
                     -t\kappa(1+\kappa)\|\zeta\|^2}
                    {2(1+\kappa)\|\txi\|^2}\right\}
 \end{align*}
 for
 \ $\txi:=(\xi_{2,3},-\xi_{1,3},\xi_{1,2})^\top\in\RR^3$ \ with \ $\txi\ne0$,
 \ where \ $\eta:=(\eta_1,\eta_2,\eta_3)^\top\in\RR^3$,
 \ $\zeta:=(\zeta_1,\zeta_2,\zeta_3)^\top\in\RR^3$, \ and
 \begin{equation}
  \kappa:=\frac{t\|\txi\|}{2}\coth\left(\frac{t\|\txi\|}{2}\right)-1,\quad
  \teta:=\frac{\sqrt{t}\kappa}{\|\txi\|^2}\xi\zeta+i\sqrt{t}\eta,\quad
   \xi:=\begin{bmatrix}
         0 & \xi_{1,2} & \xi_{1,3} \\
         -\xi_{1,2} & 0 & \xi_{2,3} \\
         -\xi_{1,3} & -\xi_{2,3} & 0
        \end{bmatrix}.
 \end{equation}
(Here \ $\|\cdot\|$ \ and \ $\langle\cdot\,,\cdot\rangle$ \ denote the
 Euclidean norm and scalar product, respectively.)
\end{Pro}

To calculate the Fourier transform of \eqref{NINE} we will use the
 constructions of the processes \ $(W^*_k(t))_{t\sgeq0}$ \ and
 \ $(W_{k,\ell}(t))_{\sgeq0}$ \ (see Lemma \ref{CHAMAU}) and the  following
 lemma.

\begin{Lem}\label{SLEM} \
Let \ $X$ \ be a \ $k$--dimensional real random vector with standard normal
 distribution.
Then we have
 \begin{equation}
  \EE\exp\big\{\langle\teta,X\rangle-s\langle BX,X\rangle\big\}
  =\frac{1}{\sqrt{\Det(I+2sB)}}
   \exp\left\{\frac{1}{2}\big\langle\teta,(I+2sB)^{-1}\teta\big\rangle\right\},
 \end{equation}
 for all \ $\teta\in\CC^k$, \ $s\in\RR^+$ \ and real symmetric positive
 semidefinite matrices \ $B$.
\ (Here \ $I$ \ denotes the \ $k\times k$ \ identity matrix.)
\end{Lem}

\noindent
\textbf{Proof.} \
Consider the decomposition \ $B=U\Lambda U^\top$, \ where \ $\Lambda$ \ is the
 \ $k\times k$ \ diagonal matrix containing the eigenvalues of \ $B$ \ in its
 diagonal and \ $U$ \ is an orthogonal matrix.
Then the random vector \ $Y:=U^\top X$ \ has also a standard
normal distribution. This implies that
 \begin{align*}
  \EE\exp\big\{\langle\teta,X\rangle-s\langle BX,X\rangle\big\}
  &=\EE\exp\big\{\langle\teta,UY\rangle -s\langle\Lambda Y,Y\rangle\big\}\\
  &=\frac{1}{\sqrt{(2\pi)^k}}\int_{\RR^k}
    \exp\left\{\langle\teta,Uy\rangle
               -s\langle \Lambda y,y\rangle
               -\frac{1}{2}\langle y,y\rangle\right\}\,\dd y,
 \end{align*}
 where \ $y=(y_1,\cdots,y_k)^\top\in\RR^k$.
\ Let \ $\lambda_1,\cdots,\lambda_k$ \ denote the eigenvalues of the matrix
 \ $B$.
\ A simple computation shows that
 \begin{align*}
  &\langle\teta,Uy\rangle-s\langle \Lambda y,y\rangle
   -\frac{1}{2}\langle y,y\rangle
   =-\sum_{j=1}^k\left(s\lambda_j+\frac{1}{2}\right)y_j^2
    +\sum_{j=1}^k(U^\top\Re\teta)_jy_j
    +i\sum_{j=1}^k(U^\top\Im \teta)_jy_j\\
  &=i\sum_{j=1}^k(U^\top\Im \teta)_jy_j
    -\sum_{j=1}^k\frac{1+2s\lambda_j}{2}
                 \left(y_j-\frac{(U^\top\Re \teta)_j}{1+2s\lambda_j}\right)^2
    +\sum_{j=1}^k\frac{(U^\top\Re\teta)_j^2}{2(1+2s\lambda_j)}.
 \end{align*}
Using the well--known formula for the Fourier transform of a standard normal
 distribution
 \begin{equation}\label{karakfg}
  \int_\RR\exp\left\{ixt-\frac{(x-m)^2}{2\sigma^2}\right\}\,\dd x
  =\sqrt{2\pi}\sigma\exp\left\{imt-\frac{1}{2}\sigma^2t^2\right\},
 \end{equation}
 for all \ $t,m\in\RR$ \ and \ $\sigma>0$, \ we obtain
 \begin{align*}
  &\EE\exp\big\{\langle\teta,X\rangle -s\langle BX,X\rangle\big\}\\
  &=\frac{1}{\sqrt{\prod_{j=1}^k(1+2s\lambda_j)}}
    \exp\left\{i\sum_{j=1}^k\frac{(U^\top\Re\teta)_j(U^\top\Im \teta)_j}
                                 {1+2s\lambda_j}
               -\sum_{j=1}^k\frac{(U^\top\Im\teta)_j^2}{2(1+2s\lambda_j)}
               +\sum_{j=1}^k\frac{(U^\top\Re\teta)_j^2}
                                 {2(1+2s\lambda_j)}\right\}.
 \end{align*}
Hence the assertion.
\proofend

\noindent
\textbf{Proof of Proposition \ref{jfEF}.} \
Because of the self--similarity property of the Wiener process, the random
 vectors
 \ $\big(W_k(t),W_{\ell}^*(t),W_{p,q}(t)
         :1\leq k,\ell\leq d, \ 1\leq p\!<\!q\leq d\big)$
 \ and
 \ $\big(c^{-1/2}W_k(ct),c^{-3/2}W_{\ell}^*(ct),c^{-1}W_{p,q}(ct)
         :1\leq k,\ell\leq d, \ 1\leq p\!<\!q\leq d\big)$
 \ have the same distribution for all \ $t\geq 0$ \ and \ $c>0$.
\ Hence
 \begin{align*}
  &\tF_t(\eta_1,\eta_2,\eta_3,\zeta_1,\zeta_2,\zeta_3,
         \xi_{1,2},\xi_{1,3},\xi_{2,3})\\
  &=\tF_{2\pi}
    \left(\sqrt{\frac{t}{2\pi}}\eta_1,\sqrt{\frac{t}{2\pi}}\eta_2,
          \sqrt{\frac{t}{2\pi}}\eta_3,
          \left(\frac{t}{2\pi}\right)^{3/2}\!\!\zeta_1,
          \left(\frac{t}{2\pi}\right)^{3/2}\!\zeta_2,
          \left(\frac{t}{2\pi}\right)^{3/2}\!\zeta_3,\frac{t}{2\pi}\xi_{1,2},
          \frac{t}{2\pi}\xi_{1,3},\frac{t}{2\pi}\xi_{2,3}\right),
 \end{align*}
 so it is sufficient to determine \ $\tF_{2\pi}$.
\ By the definition of the Fourier transform we get
 \begin{equation}\label{jfEFSEG}
  \begin{split}
   &\tF_{2\pi}(\eta_1,\eta_2,\eta_3,\zeta_1,\zeta_2,\zeta_3,
               \xi_{1,2},\xi_{1,3},\xi_{2,3})\\
   &=\EE\exp\bigg\{i\bigg(\sum_{j=1}^3\eta_jW_j(2\pi)
                          +\sum_{j=1}^3\zeta_jW_j^*(2\pi)
                          +\sum_{1\sleq j<k\sleq 3}
                            \xi_{j,k}W_{j,k}(2\pi)\bigg)\bigg\}.
  \end{split}
 \end{equation}
For abbreviation let \ $\tF_{2\pi}$ \ denote
 \ $\tF_{2\pi}(\eta_1,\eta_2,\eta_3,\zeta_1,\zeta_2,\zeta_3,
               \xi_{1,2},\xi_{1,3},\xi_{2,3})$.
\ Define the random vector \ $\chi:=(\chi_1,\chi_2,\chi_3)^\top$ \ by
 \begin{align*}
  \chi_1&:=-\xi_{1,2}\frac{1}{\sqrt{\pi}}W_2(2\pi)
           -\xi_{1,3}\frac{1}{\sqrt{\pi}}W_3(2\pi)-2\sqrt{\pi}\zeta_1,\\
  \chi_2&:=\xi_{1,2}\frac{1}{\sqrt{\pi}}W_1(2\pi)
           -\xi_{2,3}\frac{1}{\sqrt{\pi}}W_3(2\pi)-2\sqrt{\pi}\zeta_2,\\
  \chi_3&:=\xi_{1,3}\frac{1}{\sqrt{\pi}}W_1(2\pi)
           +\xi_{2,3}\frac{1}{\sqrt{\pi}}W_2(2\pi)-2\sqrt{\pi}\zeta_3.
 \end{align*}
Substituting the expressions \eqref{WCON1}, \eqref{WCON2} for
 \ $W_{j,k}(2\pi)$ \ and \ $W_\ell^*(2\pi)$ \ into the formula \eqref{jfEFSEG},
 taking conditional expectation with respect to
 \ $\{W_j(2\pi),\,a_n^{(j)},\,1\leq j\leq 3,\,n\geq 1\}$, \ and using the
 identity \ $\EE(\EE(X\vert Y))=\EE X$ \ (where \ $X$, $Y$ \ random variables,
 \ $\EE\vert X\vert<\infty$), \ we obtain
 \begin{align*}
  \tF_{2\pi}=\EE\bigg[
  &\exp\Big\{i\big(\eta_1W_1(2\pi)+\eta_2W_2(2\pi)+\eta_3W_3(2\pi)\big)\Big\}\\
  &\times\EE\bigg(\exp\bigg\{i\sum_{n=1}^\infty
                               \frac{1}{n}
                               \langle\xi\cdot a_n+\chi,b_n\rangle\bigg\}
   \,\bigg\vert\,W_j(2\pi),\,a_n^{(j)},\,1\leq j\leq 3,\,n\geq1\bigg)\bigg],
 \end{align*}
 where \ $a_n:=(a_n^{(1)},a_n^{(2)},a_n^{(3)})^\top$ \ and
 \ $b_n:=(b_n^{(1)},b_n^{(2)},b_n^{(3)})^\top$.
\ Taking into account that \ $b_n^{(1)},\,b_n^{(2)},\,b_n^{(3)}$ \ are
 independent of the condition above and of each other for all \ $n\in\NN$,
 \ using the dominated convergence theorem and the explicit formula for the
 Fourier transform of a standard normal distribution we get
 $$
  \tF_{2\pi}
  =\EE\left[\exp\Big\{i\big(\eta_1W_1(2\pi)+\eta_2W_2(2\pi)
                            +\eta_3W_3(2\pi)\big)\Big\}
            \prod_{n=1}^\infty
             \exp\bigg\{-\frac{1}{2n^2}\|\xi\cdot a_n+\chi\|^2\bigg\}\right].
 $$
Since \ $\xi$ \ is a skew symmetric matrix, there exists an orthogonal matrix
 \ $M=(m_{j,k})_{1\sleq j,k\sleq3}$ \ such that
 $$
  M^\top\xi M=\begin{bmatrix}
                 0 & p & 0 \\
                -p & 0 & 0 \\
                 0 & 0 & 0
               \end{bmatrix}=:P.
 $$
The orthogonality of \ $M$ \ implies \ $M^{-1}=M^\top$, \ hence \ $\xi M=MP$.
\ We have
 $$
  MP=\begin{bmatrix}
      -pm_{1,2} & pm_{1,1} & 0 \\
      -pm_{2,2} & pm_{2,1} & 0 \\
      -pm_{3,2} & pm_{3,1} & 0
     \end{bmatrix}
    =[ -p{\mathbf{m}}_2, p{\mathbf{m}}_1, 0 ],
 $$
where \ ${\mathbf{m}}_i$, \ $i=1,2,3$, \ denotes the column vectors of \ $M$,
 \ that is, \ $M=[{\mathbf{m}}_1,{\mathbf{m}}_2,{\mathbf{m}}_3]$.
\ Obviously, \ $\xi M=[\xi{\mathbf{m}}_1,\xi{\mathbf{m}}_2,\xi{\mathbf{m}}_3]$,
 \ hence \ $\xi\mathbf{m}_1=-p\mathbf{m}_2$,
 \ $\xi{\mathbf{m}}_2=p\mathbf{m}_1$, \ $\xi{\mathbf{m}}_3=0$.
Taking into account that \ $M$ \ is orthogonal, we have
 \ $\|{\mathbf{m}}_3\|=1$, \ hence
 $$
  {\mathbf{m}}_3=\pm\frac{1}{\sqrt{\xi_{1,2}^2+\xi_{1,3}^2+\xi_{2,3}^2}}
                    (\xi_{2,3},-\xi_{1,3},\xi_{1,2})^\top.
 $$
Moreover,
 \ $\xi^2{\mathbf{m}}_1
    =\xi(\xi{\mathbf{m}}_1)=\xi(-p{\mathbf{m}}_2)=-p^2{\mathbf{m}}_1$.
\ The only nonzero eigenvalue of  \ $\xi^2$  \ is
 \ $-(\xi_{1,2}^2+\xi_{1,3}^2+\xi_{2,3}^2)$, \ hence
\ $p=\pm\sqrt{\xi_{1,2}^2+\xi_{1,3}^2+\xi_{2,3}^2}$, \ and \ $M$ \ can be
 chosen such that \ ${\mathbf{m}}_3=\txi/\|\txi\|$, \ $p=\|\txi\|$, \ and thus
 \begin{equation}
  \label{MSPEC1}
  \langle{\mathbf{m}}_1,u\rangle^2+\langle{\mathbf{m}}_2,u\rangle^2
  =\|M^\top u\|^2-\langle{\mathbf{m}}_3,u\rangle^2
  =\|u\|^2-\frac{1}{\|\txi\|^2}\langle\txi,u\rangle^2,
 \end{equation}
 for all \ $u\in\RR^3$.
\ We also get
 $$
  -\xi^2=M\begin{bmatrix}
           \|\txi\|^2 & 0 & 0 \\
               0      & \|\txi\|^2   & 0 \\
               0    & 0  & 0
          \end{bmatrix}M^\top=:M\Lambda M^\top.
 $$
To continue the calculation of the Fourier transform of \eqref{NINE} we take
 conditional expectation with respect to \ $\{W_1(2\pi),W_2(2\pi),W_3(2\pi)\}$.
\ A special case of Lemma \ref{SLEM} \ is that
 $$
  \EE\exp\bigg\{-s\sum_{j=1}^nY_j^2\bigg\}
  =\frac{1}{\sqrt{\Det(I+2sD)}}
   \exp\bigg\{\Big\langle\big(2s^2D^{1/2}(I+2sD)^{-1}D^{1/2}-sI\big)m,
                         m\Big\rangle\bigg\}
 $$
 for all \ $s\in\RR^+$, \ where \ $Y=(Y_1,\cdots,Y_k)^\top$ \ is a
 \ $k$--dimensional random variable with normal distribution such that
 \ $\EE Y=m$ \ and \ $\var Y=D$.
\ Applying this formula for \ $Y=\xi\cdot a_n+\chi$ \ with \ $s=(2n^2)^{-1}$,
 \ $m=\chi$ \ and \ $D=\xi\cdot\xi^\top=-\xi^2=M\Lambda M^\top$ \ we get
 \begin{align*}
  \tF_{2\pi}=\EE\bigg[
    &\exp\Big\{i\big(\eta_1W_1(2\pi)+\eta_2W_2(2\pi)
                             +\eta_3W_3(2\pi)\big)\Big\}\\
   &\times\prod_{n=1}^\infty\frac{1}{\sqrt{\Det(I+n^{-2}\Lambda)}}
    \exp\bigg\{\frac{1}{2}
               \Big\langle\big(n^{-4}\sqrt{\Lambda}
                          (I+n^{-2}\Lambda)^{-1}\sqrt{\Lambda}
                          -n^{-2}I\big)M^{-1}\chi,M^{-1}\chi\Big\rangle\bigg\}\bigg].
 \end{align*}
Clearly \ $\Det(I+n^{-2}\Lambda)=(1+n^{-2}\|\txi\|^2)^2$.
\ Using that
 \ $\prod_{k=1}^\infty\frac{k^2\pi^2}{k^2\pi^2+x^2}=\frac{x}{\sinh x}$,
 \ $x\coth x -1=x^2\sum_{k=1}^\infty\frac{2}{k^2\pi^2+x^2}$, \ $x\in\RR$ \ (see
 \cite{GRA65}, formulas 1.431 and 1.421), the identity \eqref{MSPEC1} and the
 fact that \ $\langle\txi,\chi\rangle^2=4\pi\langle\zeta,\txi\rangle^2$ \ we
 obtain
 \begin{align*}
  \tF_{2\pi}=
   &\frac{\pi\|\txi\|}{\sinh(\pi\|\txi\|)}
    \exp\bigg\{-\frac{\pi^3}{\|\txi\|^2}
                \left(\frac{1}{3}-\frac{\kappa}{\pi^2\|\txi\|^2}\right)
                \langle \zeta,\txi\rangle^2\bigg\}\\
  &\times\EE\exp\bigg\{i\big(\eta_1W_1(2\pi)+\eta_2W_2(2\pi)
                             +\eta_3W_3(2\pi)\big)
                       -\frac{\kappa}{4\|\txi\|^2}\|\chi\|^2\bigg\},
 \end{align*}
 where \ $\kappa=\pi\|\txi\|\coth(\pi\|\txi\|)-1$.
\ A simple computation shows that
 \begin{align*}
  \|\chi\|^2=
   &\frac{1}{\pi}
    \Big((\xi_{1,2}^2+\xi_{1,3}^2)W_1^2(2\pi)
          +(\xi_{1,2}^2+\xi_{2,3}^2)W_2^2(2\pi)
          +(\xi_{2,3}^2+\xi_{1,3}^2)W_3^2(2\pi)\Big)
    +4\pi\|\zeta\|^2\\
   &+\frac{2}{\pi}
     \Big(\xi_{1,3}\xi_{2,3}W_1(2\pi)W_2(2\pi)
          -\xi_{1,2}\xi_{2,3}W_1(2\pi)W_3(2\pi))
          +\xi_{1,2}\xi_{1,3}W_2(2\pi)W_3(2\pi)\Big)\\
   &-4(\xi_{1,2}\zeta_2+\xi_{1,3}\zeta_3)W_1(2\pi)
    +4(\xi_{1,2}\zeta_1-\xi_{2,3}\zeta_3)W_2(2\pi)
    +4(\xi_{1,3}\zeta_1+\xi_{2,3}\zeta_2)W_3(2\pi).
 \end{align*}
Using Lemma \ref{SLEM} with
 \ $\teta=\frac{\sqrt{2\pi}\kappa}{\|\txi\|^2}\xi\zeta+i\sqrt{2\pi}\eta$,
 \ $B:=-2\xi^2$, \ $s=\frac{\kappa}{4\|\txi\|^2}$ \ and taking into account
 that \ $\sqrt{\Det(I+2sB)}=1+\kappa$ \ we conclude
 \begin{align*}
  \tF_{2\pi}=
   &\frac{\pi\|\txi\|}{(1+\kappa)\sinh(\pi\|\txi\|)}
    \exp\bigg\{-\frac{\pi^3}{\|\txi\|^2}
                \bigg(\frac{1}{3}-\frac{\kappa}{\pi^2\|\txi\|^2}\bigg)
                \langle \zeta,\txi\rangle^2\bigg\}\\
  &\times
   \exp\bigg\{-\frac{\pi\kappa}{\|\txi\|^2}\|\zeta\|^2
              +\frac{1}{2}
               \bigg\langle\teta,
                           \Big(I-\frac{\kappa}{\|\txi\|^2}\xi^2\Big)^{-1}
                           \teta\bigg\rangle\bigg\}.
 \end{align*}
Using \eqref{MSPEC1} we get
 $$
  \bigg\langle\teta,
              \Big(I-\frac{\kappa}{\|\txi\|^2}\xi^2\Big)^{-1}
              \teta\bigg\rangle
  =\frac{1}{1+\kappa}\|\teta\|^2
   +\frac{\kappa}{1+\kappa}\frac{\langle\txi,\teta\rangle^2}{\|\txi\|^2}.
 $$
Hence the assertion.
\proofend

\noindent
\textbf{Proof of Theorem \ref{SREP}.} \
We prove only the case \ $\rank(B)=3$.
\ The cases \ $\rank(B)=1$ \ and \ $\rank(B)=2$ \ can be handled in a similar
 way.
In case \ $\rank(B)=3$ \ the measure \ $\mu$ \ is absolutely continuous and so
 Proposition \ref{ACP} implies that the partial Euclidean Fourier transform
 \ $\widetilde{f}_{2,3}$ \ of the measure \ $\mu$ \ has to be calculated in
 order to obtain the Fourier transform \ $\hmu(\pi_{\pm\lambda})$.
\ Let \ $(\mu_t)_{t\sgeq0}$ \ be a Gaussian semigroup such that \ $\mu_1=\mu$
 \ and let
 \ $\rho_1:=\sigma_{1,1}\sigma_{2,2}-\sigma_{1,2}\sigma_{2,1}$,
 \ $\rho_2:=\sigma_{1,1}\sigma_{2,3}-\sigma_{1,3}\sigma_{2,1}$,
 \ $\rho_3:=\sigma_{1,2}\sigma_{2,3}-\sigma_{1,3}\sigma_{2,2}$ \ by definition.
In case \ $\rank(B)=3$, \ the representation of \ $(\mu_t)_{t\sgeq0}$ \ by the
 process \ $(Z(t))_{t\sgeq0}$ \ (see Section \ref{PREL}) gives us
 $$
  Z_1(1)=a_1+\sum_{k=1}^3\sigma_{1,k}W_k(1),\qquad
  Z_2(1)=a_2+\sum_{k=1}^3\sigma_{2,k}W_k(1),
 $$
 $$
  Z_3(1)=a_3+\sum_{k=1}^3\sigma_{3,k}W_k(1)
         +\sum_{k=1}^3(a_2\sigma_{1,k}-a_1\sigma_{2,k})W^*_k(1)
         +\rho_1W_{1,2}(1)+\rho_2W_{1,3}(1)+\rho_3W_{2,3}(1).
 $$
This implies that the (full) Euclidean Fourier transform of the measure \ $\mu$
 \ is
 \begin{align*}
  \tf(\ts_1,\ts_2,\ts_3)
   =&\EE\exp\Big\{i\big(\ts_1Z_1(1)+\ts_2Z_2(1)+\ts_3Z_3(1)\big)\Big\}
    =\exp\big\{i(\ts_1a_1+\ts_2a_2+\ts_3a_3)\big\}\\
   &\times
    \EE\exp\bigg\{i\bigg(\sum_{k=1}^{3}
                          (\sigma_{1,k}\ts_1+\sigma_{2,k}\ts_2
                           +\sigma_{3,k}\ts_3)W_k(1)
                         +\sum_{k=1}^{3}
                           (a_2\sigma_{1,k}-a_1\sigma_{2,k})\ts_3W_k^*(1)\\
   &\phantom{\times\EE\bigg[\exp\bigg\{i\bigg)}
    +\ts_3\rho_1W_{1,2}(1)+\ts_3\rho_2W_{1,3}(1)
    +\ts_3\rho_3W_{2,3}(1)\bigg)\bigg\}.
\end{align*}
Proposition \ref{ACP} shows that we may suppose \ $\ts_3\ne 0$.
\ Using Proposition \ref{jfEF} and the facts that
 \begin{equation}
  \label{SREPSEG1}
  \begin{split}
   &\sum_{k=1}^{d}(a_2\sigma_{1,k}-a_1\sigma_{2,k})^2
    =b_{2,2}a_1^2-2b_{1,2}a_1a_2+b_{1,1}a_2^2,\qquad d=1,2,3,\\
   &\rho_1(a_1\sigma_{2,3}-a_2\sigma_{1,3})
    -\rho_2(a_1\sigma_{2,2}-a_2\sigma_{1,2})
    +\rho_3(a_1\sigma_{2,1}-a_2\sigma_{1,1})=0,\\
   &\delta^2=\rho_1^2+\rho_2^2+\rho_3^2,
  \end{split}
 \end{equation}
 we get
 \begin{align*}
  \tf(\ts_1,\ts_2,\ts_3)
  =\frac{1}{\cosh(\vert\ts_3\vert\delta/2)}\exp\bigg\{
   &i\big(\ts_1a_1+\ts_2a_2+\ts_3a_3\big)
    -\frac{\kappa}{2\delta^2}
     \big(b_{2,2}a_1^2-2b_{1,2}a_1a_2+b_{1,1}a_2^2\big)\\
   &+\frac{1}{2(1+\kappa)}\|\teta\|^2
    +\frac{\kappa}{2(1+\kappa)}
     \frac{\langle\txi,\teta\rangle^2}{\delta^2}\bigg\},
 \end{align*}
 where
 $$
  \kappa=\frac{|\ts_3|\delta}{2}\coth\left(\frac{|\ts_3|\delta}{2}\right)-1,
   \qquad
  \teta=-\frac{\kappa}{\delta^2}(v_1,v_2,v_3)^\top+i\Sigma^\top\ts
 $$
 with
 \begin{align*}
  v_1&=\rho_1(a_1\sigma_{2,2}-a_2\sigma_{1,2})
       +\rho_2(a_1\sigma_{2,3}-a_2\sigma_{1,3}),\\
  v_2&=-\rho_1(a_1\sigma_{2,1}-a_2\sigma_{1,1})
       +\rho_3(a_1\sigma_{2,3}-a_2\sigma_{1,3}),\\
  v_3&=-\rho_2(a_1\sigma_{2,1}-a_2\sigma_{1,1})
       -\rho_3(a_1\sigma_{2,2}-a_2\sigma_{1,2}),
 \end{align*}
 and \ $\ts:=(\ts_1,\ts_2,\ts_3)^\top$, \ $\txi:=(\rho_3,-\rho_2,\rho_1)^\top$.
\ It can be easily checked that
 \begin{align*}
  \langle\txi,\teta\rangle^2
   &=-\ts_3^2\Det B,\\
  \|\teta\|^2
   &=-\langle B\ts,\ts\rangle+\frac{\kappa^2}{\delta^4}\langle v,v\rangle
     -2i\frac{\kappa}{\delta^2}
      \big((\ts_1a_1+\ts_2a_2)\delta^2+\ts_3(a_1\delta_3+a_2\delta_1)\big),\\
  \ts^\top B\ts
   &=b_{1,1}\left(\ts_1+\frac{b_{1,2}\ts_2+b_{1,3}\ts_3}{b_{1,1}}\right)^2
     +\frac{1}{b_{1,1}}
      \begin{bmatrix}
       \ts_2 \\ \ts_3
      \end{bmatrix}^\top
      \begin{bmatrix}
       \delta^2 & \delta_1\\
       \delta_1 & \delta_4
      \end{bmatrix}
      \begin{bmatrix}
       \ts_2 \\ \ts_3
      \end{bmatrix},
 \end{align*}
 where \ $\delta_3:=b_{1,3}b_{2,2}-b_{1,2}b_{2,3}$ \ and
 \ $\delta_4:=b_{1,1}b_{3,3}-b_{1,3}^2$.
\ Using \eqref{partialfourier}, the identities above and \eqref{karakfg}, the
 partial Fourier transform \ $\tf_{2,3}$ \ can be calculated as follows
 \begin{align*}
  &\tf_{2,3}(s_1,\ts_2,\ts_3)
   =\sqrt{\frac{\vert\ts_3\vert\delta}{2\pi b_{1,1}\sinh(|\ts_3|\delta)}}
    \exp\bigg\{-\frac{\kappa}{2(1+\kappa)\delta^2}
                (b_{2,2}a_1^2-2b_{1,2}a_1a_2+b_{1,1}a_2^2)\\
  &-\frac{\kappa}{2(1+\kappa)\delta^2}\ts_3^2\Det B
   -\frac{1}{2(1+\kappa)b_{1,1}}
    \begin{bmatrix}
     \ts_2 \\ \ts_3
    \end{bmatrix}^\top
    \begin{bmatrix}
     \delta^2 & \delta_1\\
     \delta_1 & \delta_4
    \end{bmatrix}
    \begin{bmatrix}
     \ts_2 \\ \ts_3
    \end{bmatrix}
   -\frac{1+\kappa}{2b_{1,1}}\left(\frac{a_1}{1+\kappa}-s_1\right)^2\\
  &-\frac{b_{1,2}\ts_2+b_{1,3}\ts_3}{b_{1,1}}
    \left(\frac{a_1}{1+\kappa}-s_1\right)
   +i\bigg(\ts_2a_2+\ts_3a_3
           -\frac{\kappa}{(1+\kappa)\delta^2}
            \big(\ts_2a_2\delta^2+\ts_3(a_1\delta_3+a_2\delta_1)\big)\bigg)
   \bigg\}.
 \end{align*}
Finally Proposition \ref{ACP} implies that the Fourier transform
 \ $\hmu(\pi_{\pm\lambda})$ \ is an integral operator on \ $L^2(\RR)$,
 $$
  [\hmu(\pi_{\pm\lambda})u](x)=\int_\RR K_{\pm\lambda}(x,y)u(y)\,\dd y,
 $$
 where \ $K_{\pm\lambda}$ \ has the form given in Theorem \ref{SREP}.
\proofend

\section{Convolution of Gaussian measures}

The convolution of two probability measures \ $\mu'$ \ and \ $\mu''$ \ on
 \ $\HH$ \ is defined by
 $$
  (\mu'*\mu'')(A):=\int_{\HH}\mu''(h^{-1}A)\,\mu'(\dd h),
 $$
 for all Borel sets \ $A$ \ in \ $\HH$.

First we give an explicit formula for the Fourier transform of a convolution of
 two Gaussian measures on \ $\HH$.

\begin{Thm}\label{CONV} \
Let \ $\mu'$ \ and \ $\mu''$ \ be Gaussian measures on \ $\HH$ \ with
 parameters \ $(a',B')$ \ and \ $(a'',B'')$, \ respectively.
Then we have
 \begin{align*}
  \left(\mu'*\mu''\right)\ssphat(\chi_{\alpha,\beta})
   =\exp\bigg\{&i\big((a_1'+a_1'')\alpha+(a_2'+a_2'')\beta\big)\\
               &-\frac{1}{2}
                 \Big((b_{1,1}'+b_{1,1}'')\alpha^2
                      +2(b_{1,2}'+b_{1,2}'')\alpha\beta
                      +(b_{2,2}'+b_{2,2}'')\beta^2\Big)\bigg\},
 \end{align*}
 \begin{align*}
  \left[\left(\mu'*\mu''\right)\ssphat(\pi_{\pm\lambda})u\right](x)
  =\begin{cases}
    L_{\pm\lambda}(x)u\big(x+\sqrt{\lambda}(a_1'+a_1'')\big)
     & \text{if \ $b_{1,1}'=b_{1,1}''=0$,} \\[2mm]
    \DS\int_\RR K_{\pm\lambda}(x,y)u(y)\,\dd y & \text{otherwise,}
   \end{cases}
 \end{align*}
 where \ $L_{\pm\lambda}(x)$ \ is given by
 \begin{align*}
  \exp\bigg\{&\pm i\Big(\lambda\big(a_3'+a_3''+(a_1'a_2'+a_1''a_2'')/2\big)
                        +\sqrt{\lambda}(a_2'+a_2'')x+\lambda a_1'a_2''\Big)\\
             &-\frac{\lambda^2}{2}
               \Big(b_{3,3}'+b_{3,3}''+a_1'b_{2,3}'+(2a_1'+a_1'')b_{2,3}''
                    +\big((a_1')^2b_{2,2}'+(a_1'')^2b_{2,2}''\big)/3
                    +a_1'(a_1'+a_1'')b_{2,2}''\Big)\\
            &-\frac{\lambda^{3/2}}{2}
              x\big(2b_{2,3}'+2b_{2,3}''+a_1'b_{2,2}'
                    +(2a_1'+a_1'')b_{2,2}''\big)
             -\frac{\lambda}{2}x^2(b_{2,2}'+b_{2,2}'')\bigg\},
 \end{align*}
 and \
  $K_{\pm\lambda}(x,y)
   :=C\exp\left\{-\frac{1}{2}\mathbf{z}^\top
   V\mathbf{z}\right\}$, \ $\mathbf{z}:=(x,y,1)^\top$, \
 with
 $$C:=\begin{cases}
       C_{\pm\lambda}(B')
        & \mbox{\rm if \ $b_{1,1}'>0$ \ and \ $b_{1,1}''=0$,}\\[2mm]
       C_{\pm\lambda}(B'')
        & \mbox{\rm if \ $b_{1,1}'=0$ \ and \ $b_{1,1}''>0$,}\\[2mm]
       C_{\pm\lambda}(B')C_{\pm\lambda}(B'')
        \sqrt{\frac{\DS2\pi}{\DS d_{2,2}'+d_{1,1}''}}
        & \mbox{\rm if \ $b_{1,1}'>0$ \ and \ $b_{1,1}''>0$}
      \end{cases}$$
 (taking the square root with positive real part) and
 $$V:=\begin{cases}
       D_{\pm\lambda}(a',B')
       +\begin{bmatrix}
         0 & 0                      & -\sqrt{\lambda}a_1''d_{1,2}' \\
         0 & \lambda b_{2,2}''      & p_{2,3} \\
         -\sqrt{\lambda}a_1''d_{1,2}'  & p_{3,2} & p_{3,3}
        \end{bmatrix}
       & \text{if \ $b_{1,1}'>0$ \ and \ $b_{1,1}''=0$,}\\[10mm]
       \begin{bmatrix}
        \lambda b_{2,2}'      & 0 & q_{1,3} \\
        0                     & 0 & \sqrt{\lambda}a_1'd_{1,2}'' \\
        q_{3,1} & \sqrt{\lambda}a_1'd_{1,2}'' & q_{3,3}
       \end{bmatrix}
       +D_{\pm\lambda}(a'',B'')
       & \text{if \ $b_{1,1}'=0$ \ and \ $b_{1,1}''>0$,}\\[10mm]
       \begin{bmatrix}
        d_{1,1}' & 0 & d_{1,3}' \\
        0 & d_{2,2}'' & d_{2,3}'' \\
        d_{3,1}' & d_{3,2}'' & d_{3,3}'+d_{3,3}''
       \end{bmatrix}
      -\frac{\DS UU^\top}{\DS d_{2,2}'+d_{1,1}''}
       & \text{if \ $b_{1,1}'>0$ \ and \ $b_{1,1}''>0$,}
      \end{cases}$$
 where \ $d_{j,k}':=d_{j,k}^{\pm\lambda}(a',B')$,
 \ $d_{j,k}'':=d_{j,k}^{\pm\lambda}(a'',B'')$ \ for \ $1\leq j,k\leq3$ \ and
 \begin{align*}
  U&:=(d_{1,2}',d_{2,1}'',d_{3,2}'+d_{3,1}'')^\top,\\
  p_{2,3}&:=p_{3,2}:=-\sqrt{\lambda}a_1''d_{2,2}'
                     +\lambda^{3/2}(2b_{2,3}''-a_1''b_{2,2}'')/2
                     \mp i\sqrt{\lambda}a_2'',\\
  p_{3,3}&:=-\sqrt{\lambda}a_1''(d_{2,3}'+d_{3,2}')
            +\lambda(a_1'')^2d_{2,2}'
            +\lambda^2\big(b_{3,3}''-a_1''b_{2,3}''+(a_1'')^2b_{2,2}''/3\big)
            \mp i\lambda(2a_3''-a_1''a_2''),\\
  q_{1,3}&:=q_{3,1}:=\sqrt{\lambda}a_1'd_{1,1}''
                     +\lambda^{3/2}(a_1'b_{2,2}'+2b_{2,3}')/2
                     \mp i\sqrt{\lambda}a_2',\\
  q_{3,3}&:=\sqrt{\lambda}a_1'(d_{1,3}''
            +d_{3,1}'')+\lambda (a_1')^2d_{1,1}''
            +\lambda^2\big(b_{3,3}'+a_1'b_{2,3}'+(a_1')^2b_{2,2}'/3\big)
            \mp i\lambda(2a_3'+a_1'a_2').
 \end{align*}
\end{Thm}

\noindent
\textbf{Proof.} \
If \ $b_{1,1}'>0$ \ and \ $b_{1,1}''>0$ \ then the assertion can be proved as
 in Pap \cite[Theorem 7.2]{PAP02}.
If \ $b_{1,1}'>0$ \ and \ $b_{1,1}''=0$ \ then by Theorem \ref{SREP}
 $$
  [\widehat{\mu'}(\pi_{\pm\lambda})u](x)
  =\int_{\RR}K_{\pm\lambda}'(x,y)u(y)\,\dd y
 $$
 with
 $$
  K_{\pm\lambda}'(x,y)
  :=C_{\pm\lambda}(B')
    \exp\left\{-\frac{1}{2}
                \mathbf{z}^\top D_{\pm\lambda}(a',B')\mathbf{z}\right\},
  \qquad \mathbf{z}=(x,y,1)^\top,
 $$
 and
 \begin{align*}
  [\widehat{\mu''}(\pi_{\pm\lambda})u](y)
  =\exp\bigg\{&\pm\frac{i\sqrt{\lambda}}{2}
                  \big(\sqrt{\lambda}(2a_3''+a_1''a_2'')+2a_2''y\big)
               -\frac{\lambda^2}{6}
                \big(3b_{3,3}''+3a_1''b_{2,3}''+(a_1'')^2b_{2,2}''\big)\\
              &-\frac{\lambda^{3/2}}{2}(2b_{2,3}''+a_1''b_{2,2}'')y
               -\frac{\lambda}{2}b_{2,2}''y^2\bigg\}
   u(y+\sqrt{\lambda}a_1'').
 \end{align*}
Clearly we have
 $$
  \left[\left(\mu'*\mu''\right)\ssphat(\pi_{\pm\lambda})u\right](x)
  =[\widehat{\mu'}(\pi_{\pm\lambda})\widehat{\mu''}(\pi_{\pm\lambda})u](x)
  =\int_\RR K_{\pm\lambda}'(x,y)[\widehat{\mu''}(\pi_{\pm\lambda})u](y)\,\dd y.
 $$
Using the formulas for \ $\widehat{\mu'}(\pi_{\pm\lambda})$ \ and
 \ $\widehat{\mu''}(\pi_{\pm\lambda})$ \ an easy calculation yields that
 \ $K_{\pm\lambda}$ \ has the form given in the theorem.
The other cases \ $b_{1,1}'=0,\, b_{1,1}''>0$ \ and \ $b_{1,1}'=b_{1,1}''=0$
 \ can be handled in the same way.
\proofend

We need two lemmas concerning the parameters of a Gaussian measure on \ $\HH$.

\begin{Lem}\label{GAUSSPAR1} \
Let us consider a Gaussian semigroup \ $(\mu_t)_{t\sgeq0}$ \ such that
 \ $\mu_1$ \ is a Gaussian measure on \ $\HH$ \ with parameters \ $(a,B)$.
\ Then we have
 $$
  a_i=\EE Z_i,\,\,i=1,2,3,\qquad b_{i,j}=\cov(Z_i,Z_j)\quad
  \text{if \ $(i,j)\ne(3,3)$,}
 $$
 and
 \begin{align*}
  b_{3,3}&=\var Z_3-\frac{1}{4}\Big(\var Z_1 \var Z_2-\cov(Z_1,Z_2)^2\Big)\\
         &\phantom{=\:}
          -\frac{1}{12}\Big(\var Z_2\,(\EE Z_1)^2
                             -2\cov(Z_1,Z_2)\,\EE Z_1\EE Z_2
                             +\var Z_1\,(\EE Z_2)^2\Big),
 \end{align*}
 where the distribution of the random vector \ $(Z_1,Z_2,Z_3)$ \ with values in
 \ $\RR^3$ \ is \ $\mu_1$.
\end{Lem}

\noindent
\textbf{Proof.} \
Let \ $Z(t):=(Z_1(t),Z_2(t),Z_3(t))$, \ $t\geq 0$ \ be given as in Section
 \ref{PREL}.
Taking the expectation of \ $Z(1)$ \ yields that \ $\EE(Z_i(1))=a_i$,
 \ $i=1,2,3$.
\ Using again the definition of \ $Z(1)$ \ and the fact that
 \ $B=\Sigma\cdot\Sigma^\top$ \ we get
 $$
  \var(Z_1(1))
  =\sum_{k=1}^{d}\sum_{\ell=1}^{d}
    \sigma_{1,k}\sigma_{1,\ell}\EE(W_k(1)W_{\ell}(1))
  =\sum_{k=1}^{d}\sigma_{1,k}^2=b_{1,1}.
 $$
Similar arguments show \ $\var(Z_2(1))=b_{2,2}$ \ and
 \ $\cov(Z_1(1),Z_2(1))=b_{1,2}$.
\ We also obtain
 \begin{align*}
  \cov(Z_1(1),Z_3(1))=\EE\bigg[
   &\sum_{i=1}^{d}\sigma_{1,i}W_i(1)
                  \bigg(\sum_{k=1}^{d}\sigma_{3,k}W_k(1)
                        +\sum_{k=1}^{d}
                          (a_2\sigma_{1,k}-a_1\sigma_{2,k})W_k^*(1)\bigg)\\
   &+\sum_{i=1}^{d}\sigma_{1,i}W_i(1)
     \sum_{1\sleq k<\ell\sleq d}
      (\sigma_{1,k}\sigma_{2,\ell}-\sigma_{1,\ell}\sigma_{2,k})W_{k,\ell}(1)
    \bigg],
 \end{align*}
 which implies that
 \begin{align*}
  \cov(Z_1(1),Z_3(1))=&\sum_{k=1}^{d}\sigma_{1,k}\sigma_{3,k}
                       +\sum_{i=1}^{d}\sum_{k=1}^{d}\sigma_{1,i}
                         (a_2\sigma_{1,k}-a_1\sigma_{2,k})\EE(W_i(1)W_k^*(1))\\
                      &+\sum_{i=1}^{d}\sum_{1\sleq k<\ell\sleq d}\sigma_{1,i}
                         (\sigma_{1,k}\sigma_{2,\ell}
                          -\sigma_{1,\ell}\sigma_{2,k})
                         \EE(W_i(1)W_{k,\ell}(1))=b_{1,3},
 \end{align*}
 since \ $W_i(1),\,\,1\leq i\leq d$ \ are independent of each other and
 \begin{equation}
  \label{GPARFORM1}
  \EE(W_i(1)W_k^*(1))=\EE(W_i(1)W_{k,\ell}(1))=0,\qquad
  1\leq i\leq d,\quad 1\leq k<\ell\leq d.
 \end{equation}
 Indeed,
 \begin{align*}
  \EE(W_i(1)W_k^*(1))
  &=\frac{1}{2}\lim_{n\to\infty}
    \EE\bigg[W_i(1)
             \sum_{j=1}^n\Big(W_k(s_{j-1}^{(n)})(s_j^{(n)}-s_{j-1}^{(n)})
                              -s_{j-1}^{(n)}
                               \big(W_k(s_j^{(n)})-W_k(s_{j-1}^{(n)})\big)\Big)
             \bigg],\\
  \EE(W_i(1)W_{k,\ell}(1))
  &=\frac{1}{2}\lim_{n\to\infty}
    \EE\bigg[W_i(1)
             \sum_{j=1}^n
              \Big(W_k(s_{j-1}^{(n)})
                   \big(W_{\ell}(s_j^{(n)})-W_{\ell}(s_{j-1}^{(n)})\big)\\
  &\phantom{=\frac{1}{2}\lim_{n\to\infty}\EE\bigg[W_i(1)\sum_{j=1}^n\Big(}
   -W_{\ell}(s_{j-1}^{(n)})\big(W_k(s_j^{(n)})-W_k(s_{j-1}^{(n)})\big)\Big)
   \bigg]
 \end{align*}
 for all \ $1\leq i\leq d$, \ $1\leq k<\ell\leq d$, \ where
 \ $\{s_j^{(n)}:j=0,\cdots,n\}$ \ denotes a partition of the interval \ $[0,1]$
 \ such that \ $\max_{1\sleq j\sleq n}(s_j^{(n)}-s_{j-1}^{(n)})$ \ tends to
 \ $0$ \ as \ $n$ \ goes to infinity.
We can obtain \ $\cov(Z_2(1),Z_3(1))=b_{2,3}$ \ in the same way.
Using again the form of \ $Z(t)$, \ \eqref{GPARFORM1} and the facts that
 $$
  \begin{array}{ll}
   \cov(W_{i,j}(1),W_{k,\ell}(1))=0
   &\text{for \ all \ $1\leq i<j\leq d$, \ $1\leq k<\ell\leq d$,
          \ $(i,j)\ne(k,\ell)$,}\\
   \cov(W_k^*(1),W_{\ell}^*(1))=0
   &\text{for \ all \ $1\leq k,\ell\leq d$, \ $k\ne \ell$},
  \end{array}
 $$
 we get
 \begin{align*}
  \var(Z_3(1))=&\sum_{k=1}^{d}\sigma_{3,k}^2
                +\sum_{k=1}^{d}
                  (a_2\sigma_{1,k}-a_1\sigma_{2,k})^2\var(W_k^*(1))\\
               &+\sum_{1\sleq k<\ell\sleq d}
                  (\sigma_{1,k}\sigma_{2,\ell}
                   -\sigma_{1,\ell}\sigma_{2,k})^2\var(W_{k,\ell}(1)).
 \end{align*}
L\'evy proved that the (Euclidean) Fourier transform of \ $W_{k,\ell}(1)$,
 \ $1\leq k<\ell\leq d$ \ (i.e., the characteristic function of \ $W_{k,\ell}$)
 is
 $$
  \EE\left(e^{itW_{k,\ell}(1)}\right)=\frac{1}{\cosh(t/2)},
  \qquad 1\leq k<\ell\leq d
 $$
 for all \ $t\in\RR$ \ (this follows also from Proposition \ref{jfEF}), so
 $$
  \var(W_{k,\ell}(1))
  =-\frac{\dd^2}{\dd\,t^2}\left(\frac{1}{\cosh(t/2)}\right)\bigg\vert_{t=0}
  =\frac{1}{4},\qquad 1\leq k<\ell\leq d.
 $$
Clearly \ $W_k^*$ \ has a normal distribution with zero mean and with variance
 \ $\var(W_k^*(1))=\frac{1}{12}$, \ $1\leq k\leq d$.
\ Using \eqref{SREPSEG1} we have
 $$
  \var(Z_3(1))
  =b_{3,3}+\frac{1}{4}(b_{1,1}b_{2,2}-b_{1,2}^2)
   +\frac{1}{12}(a_1^2b_{2,2}-2a_1a_2b_{1,2}+a_2^2b_{1,1}).
 $$
Hence the assertion.
\proofend

\begin{Lem}\label{GAUSSPAR2} \
Let \ $\mu'$ \ and \ $\mu''$ \ be Gaussian measures on \ $\HH$ \ with
 parameters \ $(a',B')$ \ and \ $(a'',B'')$, \ respectively.
If the convolution \ $\mu'*\mu''$ \ is a Gaussian measure on \ $\HH$ \ with
 parameters \ $(a,B)$ \ then we have
 \begin{align*}
  a_1&=a_1'+a_1'',\qquad a_2=a_2'+a_2'',\qquad
  a_3=a_3'+a_3''+\frac{1}{2}(a_1'a_2''-a_2'a_1''),\\
  b_{1,1}&=b_{1,1}'+b_{1,1}'',\qquad b_{1,2}=b_{1,2}'+b_{1,2}'',\qquad
  b_{2,2}=b_{2,2}'+b_{2,2}'',\\
  b_{1,3}&=b_{1,3}'+b_{1,3}''
           +\frac{1}{2}\big(a_2''b_{1,1}'-a_1''b_{1,2}'
                            +a_1'b_{1,2}''-a_2'b_{1,1}''\big),\\
  b_{2,3}&=b_{2,3}'+b_{2,3}''
           +\frac{1}{2}\big(a_2''b_{1,2}'-a_1''b_{2,2}'
                            +a_1'b_{2,2}''-a_2'b_{1,2}''\big),\\
  b_{3,3}&=b_{3,3}'+b_{3,3}''+a_2''b_{1,3}'-a_1''b_{2,3}'+a_1'b_{2,3}''
           -a_2'b_{1,3}''\\
         &\phantom{=\:}
          +\frac{1}{6}
           \Big(-a_1'a_1''b_{2,2}'+(a_1'')^2b_{2,2}'+(a_1')^2b_{2,2}''
                -a_1'a_2''b_{2,2}''+a_1'a_2''b_{1,2}'+a_1''a_2'b_{1,2}'
                -2a_1''a_2''b_{1,2}'\\
         &\phantom{=\:+\frac{1}{6}\Big(}
           -2a_1'a_2'b_{1,2}''+a_1'a_2''b_{1,2}''+a_1''a_2'b_{1,2}''
           -a_2'a_2''b_{1,1}'+(a_2'')^2b_{1,1}'+(a_2')^2b_{1,1}''
           -a_2'a_2''b_{1,1}''\Big).
 \end{align*}
\end{Lem}

\noindent
\textbf{Proof.} \
Let \ $Z'=(Z'_1,Z'_2,Z'_3)^\top$ \ and \ $Z''=(Z''_1,Z''_2,Z''_3)^\top$ \ be
 independent random variables with values in \ $\RR^3$ \ such that the
 distribution of \ $Z'$ \ is \ $\mu'$ \ and the distribution of \ $Z''$ \ is
 \ $\mu''$, \ respectively.
Then the convolution \ $\mu'*\mu''$ \ is the distribution of the random
 variable
 $$
  \Big(Z'_1+Z''_1,Z'_2+Z_2'',
       Z'_3+Z_3''+\frac{1}{2}(Z_1'Z_2''-Z_1''Z_2')\Big)=:(Z_1,Z_2,Z_3).
 $$
Using Lemma \ref{GAUSSPAR1} we get
 \begin{align*}
  &a_1=\EE Z_1=\EE Z_1'+\EE Z_1''=a_1'+a_1'',\\
  &a_2=\EE Z_2=\EE Z_2'+\EE Z_2''=a_2'+a_2'',\\
  &a_3=\EE Z_3=\EE Z_3'+\EE Z_3''
       +\frac{1}{2}\big(\EE Z_1'\EE Z_2''-\EE Z_1'' \EE Z_2'\big)
      =a_3'+a_3''+\frac{1}{2}(a_1'a_2''-a_2'a_1''),
 \end{align*}
 since \ $Z'$ \ and \ $Z''$ \ are independent of each other.
Similar arguments show that
 \begin{align*}
  &b_{1,1}=\var Z_1=\var Z_1'+\var Z_1''=b_{1,1}'+b_{1,1}'',\\
  &b_{2,2}=\var Z_2=\var Z_2'+\var Z_2''=b_{2,2}'+b_{2,2}'',\\
  &b_{1,2}=\cov(Z_1,Z_2)=b_{1,2}'+b_{1,2}''.
 \end{align*}
We also have
 \begin{align*}
  b_{1,3}=\cov(Z_1,Z_3)
  &=\cov(Z'_1,Z'_3)+\cov(Z''_1,Z''_3)\\
  &+\frac{1}{2}\Big(\cov(Z_1',Z_1'Z_2'')-\cov(Z_1',Z_2'Z_1'')
                    +\cov(Z_1'',Z_1'Z_2'')-\cov(Z_1'',Z_1''Z_2')\Big).
 \end{align*}
Using this and Lemma \ref{GAUSSPAR1} the validity of the formula for
 \ $b_{1,3}$ \ can be easily checked.
For example we have
 $$
  \cov(Z'_1,Z'_1Z''_2)=\EE\big((Z_1')^2Z_2''\big)-\EE Z_1'\EE(Z_1'Z_2'')
  =\big(b_{1,1}'+(a_1')^2\big)a_2''-(a_1')^2a_2''=a_2''b_{1,1}'.
 $$
The validity of the formula for \ $b_{2,3}$ \ can be proved in the same way.
Lemma \ref{GAUSSPAR1} implies that
 \begin{align*}
  \var Z_3
  &=b_{3,3}+\frac{1}{4}(b_{1,1}b_{2,2}-b_{1,2}^2)
    +\frac{1}{12}(a_1^2b_{2,2}-2a_1a_2b_{1,2}+a_2^2b_{1,1})
   =\cov(Z_3,Z_3)\\
  &=\cov(Z_3',Z_3')+\cov(Z_3'',Z_3'')+\cov(Z_3',Z_1'Z_2'')
    -\cov(Z_3',Z_1''Z_2')+\cov(Z_3'',Z_1'Z_2'')\\
  &\phantom{=\:}
    -\cov(Z_3'',Z_1''Z_2')
    +\frac{1}{4}
     \Big(\cov(Z_1'Z_2'',Z_1'Z_2'')-\cov(Z_1'Z_2'',Z_1''Z_2')\\
  &\phantom{=\:-\cov(Z_3'',Z_1''Z_2')+\frac{1}{4}\bigg(}
    -\cov(Z_1''Z_2',Z_1'Z_2'')+\cov(Z_1''Z_2',Z_1''Z_2')\Big).
 \end{align*}
Using again Lemma \ref{GAUSSPAR1} and substituting the formulas for
 \ $b_{1,1}$, $b_{1,2}$, $b_{2,2}$, $a_1$ \ and \ $a_2$ \ into the formula
 above, an easy calculation shows the validity of the formula for \ $b_{3,3}$.
\proofend

Our aim is to give necessary and sufficient conditions for a convolution of two
 Gaussian measures to be a Gaussian measure.
Using the fact that the Fourier transform is injective (i.e., if \ $\mu$ \ and
 \ $\nu$ \ are probability measures on \ $\HH$ \ such that
 \ $\hmu(\chi_{\alpha,\beta})=\widehat{\nu}(\chi_{\alpha,\beta})$
 \ for all \ $\alpha,\beta\in\RR$ \ and
 \ $\hmu(\pi_{\pm\lambda})=\widehat{\nu}(\pi_{\pm\lambda})$ \ for all
 \ $\lambda>0$ \ then \ $\mu=\nu$), our task can be fulfilled in the following
 way.
We take the Fourier transform of the convolution of two Gaussian measures
 \ $\mu'$ \ and \ $\mu''$ \ with parameters \ $(a',B')$ \ and \ $(a'',B'')$
 \ at all one--dimensional and at all Schr\"odinger representations and then we
 search for necessary and sufficient conditions under which this Fourier
 transform has the form given in Theorem \ref{SREP}.
First we sketch our approach to obtain necessary conditions.
By Theorem \ref{CONV}, \ $\left(\mu'*\mu''\right)\ssphat(\pi_{\pm\lambda})$
 \ is an integral operator for \ $b_{1,1}'+b_{1,1}''>0$, \ and it is a product
 of certain shift and multiplication operators for \ $b_{1,1}'+b_{1,1}''=0$.
\ If the convolution \ $\mu'*\mu''$ \ is a Gaussian measure with parameters
 \ $(a,B)$ \ then, by Theorem \ref{SREP},
 \ $\left(\mu'*\mu''\right)\ssphat(\pi_{\pm\lambda})$ \ is an integral operator
 for \ $b_{1,1}>0$, \ and it is a product of certain shift and multiplication
 operators for \ $b_{1,1}=0$.
\ By Lemma \ref{GAUSSPAR2}, we have \ $b_{1,1}=b_{1,1}'+b_{1,1}''$, \ hence
 \ $b_{1,1}=0$ \ if and only if \ $b_{1,1}'+b_{1,1}''=0$. \
 Hence if \ $b_{1,1}>0,$ \ the integral operator
 \ $\left(\mu'*\mu''\right)\ssphat(\pi_{\pm\lambda})$ \ can be
 written with the kernel function given in Theorem \ref{SREP} and
 also with the kernel function given in Theorem \ref{CONV}.
 In the next lemma we derive some consequences of this observation.

\begin{Lem}\label{LEM_INT_OP}
 Let \ $\mu'$ \ and \ $\mu''$ \ be Gaussian measures on \ $\HH$ \ with
 parameters \ $(a',B')$ \ and \ $(a'',B''),$ \ respectively.
 Suppose that \ $\mu'*\mu''$ \ is a Gaussian measure on \ $\HH$ \
 with parameters \ $a=(a_i)_{1\sleq i\sleq 3},$ \
 $B=(b_{j,k})_{1\sleq j,k\sleq 3}$ \ such that \ $b_{1,1}>0.$ \ Then
 \ $d_{j,k}^{\pm\lambda}=v_{j,k}^{\pm\lambda}$ \ for all
  \ $1\leq j,k\leq 3$ \ with \ $(j,k)\ne(3,3)$ \ and for all \ $\lambda>0$, \
and
$$
    C_{\pm\lambda}(B)\exp\left\{-\frac{1}{2}d_{3,3}^{\pm\lambda}\right\}
       =C\exp\left\{-\frac{1}{2}v_{3,3}^{\pm\lambda}\right\},
     \quad \lambda>0,
$$
 where \ $C_{\pm\lambda}(B),$ $d_{j,k}^{\pm\lambda}:=d_{j,k}^{\pm\lambda}(a,B)$,
 \ $1\leq j,k\leq 3$ \ and
 \ $C,$ $V=:(v_{j,k}^{\pm\lambda})_{1\sleq j,k\sleq 3}$ \ are defined
 in Theorems \ref{SREP} and \ref{CONV}, respectively.
\end{Lem}

\noindent{\bf Proof.}
 The Fourier transform \ $\left(\mu'*\mu''\right)\ssphat(\pi_{\pm\lambda})$ \
 is a bounded linear operator on \ $L^2(\RR),$ \ and since \ $b_{1,1}>0,$
 \ Theorem \ref{SREP} yields that it is an integral operator on \ $L^2(\RR),$ \
\begin{align}\label{LEM_INT_OP1}
  \left[\left(\mu'*\mu''\right)\ssphat(\pi_{\pm\lambda})u\right](x)
      =\int_\RR K_{\pm\lambda}(x,y)u(y)\;\dd y,\quad
      u\in L^2(\RR),\;\;x\in\RR,
\end{align}
where
$$
  K_{\pm\lambda}(x,y)
     =C_{\pm\lambda}(B)\exp\left\{-\frac{1}{2}\mathbf{z}^\top
                    D_{\pm\lambda}(a,B)\mathbf{z}\right\},
    \quad \mathbf{z}=(x,y,1)^\top.
$$
Let us write
 \ $d_{j,k}'=:d_{j,k}^{\pm\lambda}(a',B')$ \ and
 \ $d_{j,k}''=:d_{j,k}^{\pm\lambda}(a'',B'')$ \
 for \ $1\leq j,k\leq 3$ \ as in Theorem \ref{CONV}.
 By Lemma \ref{GAUSSPAR2}, we have \ $b_{1,1}=b_{1,1}'+b_{1,1}''$, \
 hence \ $b_{1,1}>0$ \ implies that \ $b_{1,1}'>0$ \ or \ $b_{1,1}''>0.$ \
 Using Theorem \ref{CONV} we have
 \begin{align}\label{LEM_INT_OP2}
  \left[\left(\mu'*\mu''\right)\ssphat(\pi_{\pm\lambda})u\right](x)
      =\int_\RR \widetilde K_{\pm\lambda}(x,y)u(y)\;\dd y,\quad
      u\in L^2(\RR),\;\;x\in\RR,
\end{align}
where
 $$
  \widetilde K_{\pm\lambda}(x,y)
     =C\exp\left\{-\frac{1}{2}\mathbf{z}^\top V\mathbf{z}\right\},
    \quad \mathbf{z}=(x,y,1)^\top.
$$
 Using \eqref{LEM_INT_OP1} and \eqref{LEM_INT_OP2}, we have
 $$
   0=\int_\RR \big(K_{\pm\lambda}(x,y)-\widetilde K_{\pm\lambda}(x,y)\big)u(y)
                   \;\dd y,\quad u\in L^2(\RR),\;\;x\in\RR.
 $$
  We show that if
 \begin{align}\label{LEM_INT_OP3}
    \int_\RR\vert K_{\pm\lambda}(x,y)\vert^2\;\dd y<\infty,
      \quad
     \int_\RR\vert\widetilde K_{\pm\lambda}(x,y)\vert^2\;\dd y<\infty,
     \quad x\in\RR,
 \end{align}
 then \ $K_{\pm\lambda}(x,y)=\widetilde K_{\pm\lambda}(x,y),$ \
 $x,y\in\RR.$ \ Indeed, for all \ $x\in\RR$, \ the function
 \ $y\in\RR\mapsto K_{\pm\lambda}(x,y)-\widetilde K_{\pm\lambda}(x,y)$ \
 is in \ $L^2(\RR).$ \ Hence
 $$
  0=\int_\RR\vert K_{\pm\lambda}(x,y)-\widetilde K_{\pm\lambda}(x,y)\vert^2
                      \;\dd y,\quad x\in\RR.
 $$
 Then we get
 $$
     \int_\RR\int_\RR\vert K_{\pm\lambda}(x,y)-\widetilde K_{\pm\lambda}(x,y)\vert^2
                      \;\dd x\,\dd y=0,
 $$
 which implies that
  \ $K_{\pm\lambda}(x,y)=\widetilde K_{\pm\lambda}(x,y)$ \
 for almost every \ $x,y\in\RR.$ \
 Using that \ $K_{\pm\lambda}$ \ and \ $\widetilde K_{\pm\lambda}$ \
 are continuous, we get \ $K_{\pm\lambda}(x,y)=\widetilde K_{\pm\lambda}(x,y),$ \
 $x,y\in\RR.$ \ Now we check that \eqref{LEM_INT_OP3} is satisfied.
 Using the forms of \ $K_{\pm\lambda}$ \ and \ $\widetilde K_{\pm\lambda},$ \
 it is enough to check that
 \begin{align}\label{LEM_INT_OP4}
   &\int_\RR\exp\left\{-\mathbf{z}^\top\Re(D_{\pm\lambda}(a,B))\mathbf{z}\right\}
        \;\dd y<\infty,\quad x\in\RR, \\\label{LEM_INT_OP5}
   &\int_\RR\exp\left\{-\mathbf{z}^\top\Re(V)\mathbf{z}\right\}
         \;\dd y<\infty,\quad x\in\RR,
 \end{align}
 where \ $\mathbf{z}=(x,y,1)^\top.$ \
 Here \ $\Re(D_{\pm\lambda}(a,B))$ \ and \ $\Re(V)$ \ are real, symmetric matrices.
  Let us consider an arbitrary real, symmetric matrix
 \ $M=(m_{i,j})_{1\sleq i,j\sleq 3}$ \ with \ $m_{2,2}>0.$ \ Then
 \begin{align*}
   \mathbf{z}^\top M\mathbf{z}
      &=m_{1,1}x^2+2m_{1,2}xy+m_{2,2}y^2+2m_{1,3}x+2m_{2,3}y+m_{3,3}\\
      &=\!\left(\sqrt{m_{2,2}}y+\frac{1}{\sqrt{m_{2,2}}}(m_{1,2}x+m_{2,3})\right)^2
         \!-\frac{1}{m_{2,2}}(m_{1,2}x+m_{2,3})^2\\
      &\phantom{=\;\;}
        +m_{1,1}x^2+2m_{1,3}x+m_{3,3}.
 \end{align*}
Hence
\begin{align*}
  &\int_\RR\exp\left\{-\mathbf{z}^\top M\mathbf{z}\right\}\;\dd y
      =\exp\left\{\frac{1}{m_{2,2}}(m_{1,2}x+m_{2,3})^2
                  -m_{1,1}x^2-2m_{1,3}x-m_{3,3}\right\}\\
  &\phantom{\int_\RR\exp\left\{-\mathbf{z}^\top M\mathbf{z}
                     \right\}\;\dd y=\;\;}\times
           \int_\RR\exp\left\{-\Big(\sqrt{m_{2,2}}y
            +\frac{1}{\sqrt{m_{2,2}}}(m_{1,2}x+m_{2,3})\Big)^2\right\}\;\dd y\\
  &\phantom{\int_\RR\exp\left\{-\mathbf{z}^\top M\mathbf{z}\right\}\;\dd y\;\;}
   =\exp\left\{\frac{1}{m_{2,2}}(m_{1,2}x+m_{2,3})^2
                  -m_{1,1}x^2-2m_{1,3}x-m_{3,3}\right\}\\
  &\phantom{\int_\RR\exp\left\{-\mathbf{z}^\top M\mathbf{z}
                     \right\}\;\dd y=\;\;}
    \times\frac{1}{\sqrt{2m_{2,2}}}\int_\RR\exp\left\{-\frac{t^2}{2}\right\}\;\dd t\\
  &=\sqrt{\frac{\pi}{m_{2,2}}}
        \exp\Bigg\{\frac{1}{m_{2,2}}(m_{1,2}x+m_{2,3})^2
               -m_{1,1}x^2-2m_{1,3}x-m_{3,3}\Bigg\},
\end{align*}
 which yields that
 $$
    \int_\RR\exp\left\{-\mathbf{z}^\top M\mathbf{z}\right\}\;\dd y
       <\infty,\quad x\in\RR.
 $$
 Hence in order to prove that \eqref{LEM_INT_OP4} and \eqref{LEM_INT_OP4}
 are valid we only have to check that the \ $(2,2)$-entries of the
 matrices \ $\Re(D_{\pm\lambda}(a,B))$ \ and \ $\Re(V)$ \
 are positive. For example, if \ $b_{1,1}'>0$ \ and
 \ $b_{1,1}''>0,$ \ then
 $$
  \big(\Re(V)\big)_{2,2}
      =\Re(d_{2,2}'')-\Re\left(\frac{(d_{2,1}'')^2}{d_{2,2}'+d_{1,1}''}\right).
 $$
  If \ $b_{1,1}'b_{2,2}'-(b_{1,2}')^2=b_{1,1}''b_{2,2}''-(b_{1,2}'')^2=0,$ \
 then
 $$
    \big(\Re(V)\big)_{2,2}
       =\frac{1}{\lambda b_{1,1}''}
         -\frac{1}{\lambda^2 (b_{1,1}'')^2}
            \frac{\frac{1}{\lambda b_{1,1}'}+\frac{1}{\lambda b_{1,1}''}}
             {\left(\frac{1}{\lambda b_{1,1}'}+\frac{1}{\lambda b_{1,1}''}\right)^2
               +\left(\frac{b_{1,2}''}{b_{1,1}''}-\frac{b_{1,2}'}{b_{1,1}'}\right)^2}.
 $$
 Hence \ $\big(\Re(V)\big)_{2,2}>0$ \ if and only if
 $$
   \lambda b_{1,1}''\left[\left(\frac{1}{\lambda b_{1,1}'}+\frac{1}{\lambda b_{1,1}''}\right)^2
       +\left(\frac{b_{1,2}''}{b_{1,1}''}-\frac{b_{1,2}'}{b_{1,1}'}\right)^2\right]
    >\frac{1}{\lambda b_{1,1}'}+\frac{1}{\lambda b_{1,1}''}.
 $$
 A simple calculation shows that the latter inequality is equivalent to
 $$
   \frac{b_{1,1}''}{b_{1,1}'}\left(\frac{1}{\lambda b_{1,1}'}
                                    +\frac{1}{\lambda b_{1,1}''}\right)
  +\lambda b_{1,1}''\left(\frac{b_{1,2}''}{b_{1,1}''}-\frac{b_{1,2}'}{b_{1,1}'}\right)^2
      >0,
 $$
 which holds since \ $b_{1,1}'>0,$ \ $b_{1,1}''>0$ \ and \ $\lambda>0$. \
 The other cases can be handled similarly. Hence \eqref{LEM_INT_OP4} and
 \eqref{LEM_INT_OP5} are satisfied, and then
  \ $K_{\pm\lambda}(x,y)=\widetilde K_{\pm\lambda}(x,y),$ \ $x,y\in\RR.$ \

 Using the forms of \ $K_{\pm\lambda}$ \ and
 \ $\widetilde K_{\pm\lambda},$ \ we get
 $$
   C_{\pm\lambda}(B)\exp\left\{-\frac{1}{2}\mathbf{z}^\top
                    D_{\pm\lambda}(a,B)\mathbf{z}\right\}
     =C\exp\left\{-\frac{1}{2}\mathbf{z}^\top V\mathbf{z}\right\},
      \quad \mathbf{z}=(x,y,1)^\top.
 $$
Putting \ $\mathbf{z}=(0,0,1)^\top$ \ gives
 \begin{align}\label{LEM_INT_OP6}
    C_{\pm\lambda}(B)\exp\left\{-\frac{1}{2}d_{3,3}^{\pm\lambda}\right\}
     =C\exp\left\{-\frac{1}{2}v_{3,3}^{\pm\lambda}\right\}.
 \end{align}
Substituting \ $\mathbf{z}=(1,0,1)^\top$ \ implies
 \begin{align*}
      C_{\pm\lambda}(B)\exp\left\{-\frac{1}{2}\big(d_{1,1}^{\pm\lambda}
            +2d_{1,3}^{\pm\lambda}+d_{3,3}^{\pm\lambda}\big)\right\}
       =C\exp\left\{-\frac{1}{2}\big(v_{1,1}^{\pm\lambda}
            +2v_{1,3}^{\pm\lambda}+v_{3,3}^{\pm\lambda}\big)\right\}.
 \end{align*}
Using \eqref{LEM_INT_OP6} we have
\begin{align}\label{LEM_INT_OP7}
   d_{1,1}^{\pm\lambda}+2d_{1,3}^{\pm\lambda}
      =v_{1,1}^{\pm\lambda}+2v_{1,3}^{\pm\lambda}.
\end{align}
With \ $\mathbf{z}=(0,1,1)^\top$ \ a similar argument shows that
\begin{align}\label{LEM_INT_OP8}
   d_{2,2}^{\pm\lambda}+2d_{2,3}^{\pm\lambda}
      =v_{2,2}^{\pm\lambda}+2v_{2,3}^{\pm\lambda}.
\end{align}
Putting \ $\mathbf{z}=(1,1,1)^\top$ \ and using \eqref{LEM_INT_OP6}
we obtain
\begin{align}\label{LEM_INT_OP9}
   \begin{split}
      d_{1,1}^{\pm\lambda}+2d_{1,2}^{\pm\lambda}
        +2d_{1,3}^{\pm\lambda}
        &+d_{2,2}^{\pm\lambda}+2d_{2,3}^{\pm\lambda}\\
      &=v_{1,1}^{\pm\lambda}+2v_{1,2}^{\pm\lambda}
         +2v_{1,3}^{\pm\lambda}+v_{2,2}^{\pm\lambda}
          +2v_{2,3}^{\pm\lambda}.
   \end{split}
\end{align}
Using \eqref{LEM_INT_OP7},\eqref{LEM_INT_OP8} and
\eqref{LEM_INT_OP9}, we have \
$d_{1,2}^{\pm\lambda}=v_{1,2}^{\pm\lambda}.$ \
 If \ $\mathbf{z}=(2,0,1)^\top$ \ then using \eqref{LEM_INT_OP6} we have
$$
  d_{1,1}^{\pm\lambda}+d_{1,3}^{\pm\lambda}
      =v_{1,1}^{\pm\lambda}+v_{1,3}^{\pm\lambda}.
$$
 Using \eqref{LEM_INT_OP7} we have \ $d_{1,3}^{\pm\lambda}=v_{1,3}^{\pm\lambda}.$ \
 If \ $\mathbf{z}=(0,2,1)^\top$ \ then
$$
  d_{2,2}^{\pm\lambda}+d_{2,3}^{\pm\lambda}
      =v_{2,2}^{\pm\lambda}+v_{2,3}^{\pm\lambda}.
$$
 Using \eqref{LEM_INT_OP8} we have \ $d_{2,3}^{\pm\lambda}=v_{2,3}^{\pm\lambda}.$
\proofend

Using Lemma \ref{LEM_INT_OP} we derive necessary conditions for a
convolution of two Gaussian measures to be a Gaussian measure and
prove that they are also sufficient. The above train of thoughts
will be used in the proof of Proposition
 \ref{CENTRALGAUSS} and Theorem \ref{SCHREPPREL}.

\begin{Rem}\label{CENTRALGAUSSREM} \
By Lemma \ref{GAUSSSUPPORT}, it can be easily checked that a Gaussian measure
 \ $\mu$ \ admits parameters \ $(a,B)$ \ with \ $b_{j,k}=0$ \ for
 \ $1\leq j,k\leq 3$ \ with \ $(j,k)\ne(3,3)$ \ and \ $a_1=a_2=0$ \ if and only
 if the support of \ $\mu$ \ is contained in the center of \ $\HH$.
\end{Rem}

Now we can derive a special case of Theorem \ref{SCHREPPREL} which will be used
 in the proof of Theorem \ref{SCHREPPREL}.

\begin{Pro}\label{CENTRALGAUSS}
If \ $\mu''$ \ is a Gaussian measure on \ $\HH$ \ with parameters \
$(a'',B'')$ \ such that the support of \ $\mu''$ \ is contained in
the center of \ $\HH$ \ then for all Gaussian
 measures \ $\mu'$, \ the convolutions \ $\mu'*\mu''$ \ and \ $\mu''*\mu'$
 \ are Gaussian measures with parameters \ $(a'+a'',B'+B'')$, \ and
 \ $\mu'*\mu''=\mu''*\mu'$.
\end{Pro}

\noindent
\textbf{Proof.} \
Let \ $\mu$ \ be a Gaussian measure with parameters \ $(a'+a'',B'+B'')$.
\ By the injectivity of the Fourier transform, in order to prove that
 \ $\mu'*\mu''=\mu$ \ is valid, it is sufficient to show that
 \ $\left(\mu'*\mu''\right)\ssphat(\chi_{\alpha,\beta})
    =\hmu(\chi_{\alpha,\beta})$
 \ for all \ $\alpha,\,\beta>0$ and
 \ $\left(\mu'*\mu''\right)\ssphat(\pi_{\pm\lambda})=\hmu(\pi_{\pm\lambda})$
 \ for all \ $\lambda>0$.
\ Theorem \ref{CONV} implies that
 \ $\left(\mu'*\mu''\right)\ssphat(\chi_{\alpha,\beta})
    =\hmu(\chi_{\alpha,\beta})$
\ is valid for all one--dimensional representations \ $\chi_{\alpha,\beta}$,
 \ $\alpha,\,\beta\in\RR$.
\ Suppose that \ $b_{1,1}'\ne 0$ \ and \ $b_{1,1}'b_{2,2}'-(b_{1,2}')^2\ne 0$.
\ By Theorem \ref{CONV}, to prove
 \ $\left(\mu'*\mu''\right)\ssphat(\pi_{\pm\lambda})=\hmu(\pi_{\pm\lambda})$
 \ for all \ $\lambda>0$ \ it is sufficient to show that
 $$
  D_{\pm\lambda}(a',B')+\begin{bmatrix}
                         0 & 0 & 0 \\
                         0 & 0 & 0 \\
                         0 & 0 & \lambda^2b_{3,3}''\mp 2i\lambda a_3''
                        \end{bmatrix}
  =D_{\pm\lambda}(a'+a'',B'+B'')
 $$
 for all \ $\lambda>0$.
\ Since \ $b_{j,k}''=0$ \ for \ $1\leq j,k\leq 3$ \ with \
$(j,k)\ne(3,3)$,
 \ we have \ $d_{j,k}^{\pm\lambda}(a'+a'',B'+B'')=d_{j,k}^{\pm\lambda}(a',B')$
 \ for \ $1\leq j,k\leq 3$ \ with \ $(j,k)\ne(3,3)$.
\ So we have to check only that
 $$
  d_{3,3}^{\pm\lambda}(a',B')+\lambda^2 b_{3,3}''\mp 2i\lambda a_3''
  =d_{3,3}^{\pm\lambda}(a'+a'',B'+B'')
 $$
 for all \ $\lambda>0$.
\ Theorem \ref{SREP} implies this.
The case \ $b_{1,1}'\ne 0$, \ $b_{1,1}'b_{2,2}'-(b_{1,2}')^2=0$ \ can be proved
 similarly.
Suppose that \ $b_{1,1}'=b_{1,1}''=0$.
\ Using again Theorem \ref{SREP}, we have
 \begin{align*}
  \left[\widehat{\mu''}(\pi_{\pm\lambda})u\right](x)
   &=\exp\bigg\{\pm i\lambda a_3''-\frac{\lambda^2}{2}b_{3,3}''\bigg\}u(x),\\
  \left[\widehat{\mu'}(\pi_{\pm\lambda})u\right](x)
   &=\exp\bigg\{\pm\frac{i\sqrt{\lambda}}{2}
                   \big(\sqrt{\lambda}(2a_3'+a'_1a'_2)+2a_2'x\big)
                -\frac{\lambda^2}{6}
                 \big(3b_{3,3}'+3a_1'b_{2,3}'+(a_1')^2b_{2,2}'\big)\\
   &\phantom{=\:\exp\bigg\{}
                -\frac{\lambda^{3/2}}{2}\big(2b_{2,3}'+a_1'b_{2,2}'\big)x
                -\frac{\lambda}{2}b_{2,2}'x^2\bigg\}u(x+\sqrt{\lambda}a_1').
 \end{align*}
Theorem \ref{SREP} implies that
 \ $\left[\hmu(\pi_{\pm\lambda})u\right](x)
    =\left[\left(\mu'*\mu''\right)\ssphat(\pi_{\pm\lambda})u\right](x)$
 \ for all \ $\lambda>0$, \ $u\in L^2(\RR)$ \ and \ $x\in\RR$.
\ Hence the assertion.
\proofend

Now we give necessary and sufficient conditions under which the convolution of
 two Gaussian measures is a Gaussian measure.

\begin{Thm}\label{SCHREPPREL} \
Let \ $\mu'$ \ and \ $\mu''$ \ be Gaussian measures on \ $\HH$ \ with
 parameters \ $a'=(a_i')_{1\sleq i\sleq 3}$,
 \ $B'=(b_{j,k}')_{1\sleq j,k\sleq 3}$ \ and \ $a''=(a_i'')_{1\sleq i\sleq 3}$,
 \ $B''=(b_{j,k}'')_{1\sleq j,k\sleq 3}$, \ respectively.
Then the convolution \ $\mu'*\mu''$ \ is a Gaussian measure on \ $\HH$ \ if and
 only if one of the following conditions holds:
\renewcommand{\labelenumi}{{\rm($\widetilde{\mbox{\rm C}}$\arabic{enumi})}}
 \begin{enumerate}
  \item $b_{1,1}'>0$, \ $\delta'>0$, \ $b_{1,1}''>0$, \ $\delta''>0$, \ and
         there exists \ $\varrho>0$ \ such that \ $b_{j,k}''=\varrho b_{j,k}'$
         \ for \ $1\leq j,k\leq3$ \ with \ $(j,k)\ne(3,3)$ \ and
         \ $a_i''=\varrho a_i'$ \ for \ $i=1,2$,
  \item $b_{1,1}'>0$, \ $\delta'=0$, \ $b_{1,1}''>0$, \ $\delta''=0$, \ and
         there exists \ $\varrho>0$ \ such that \ $b_{j,k}''=\varrho b_{j,k}'$
         \ for \ $1\leq j,k\leq2$,
  \item $b_{1,1}'>0$, \ $\delta'>0$, \ $b_{j,k}''=0$ \ for \ $1\leq j,k\leq3$
        \ with \ $(j,k)\ne(3,3)$ \ and \ $a_i''=0$ \ for \ $i=1,2$,
  \item $b_{1,1}'>0$, \ $\delta'=0$, \ $b_{j,k}''=0$ \ for \ $1\leq j,k\leq3$
         \ with \ $(j,k)\ne(3,3)$,
  \item $b_{1,1}''>0$, \ $\delta''>0$, \ $b_{j,k}'=0$ \ for \ $1\leq j,k\leq3$
        \ with \ $(j,k)\ne(3,3)$ \ and \ $a_i'=0$ \ for \ $i=1,2$,
  \item $b_{1,1}''>0$, \ $\delta''=0$, \ $b_{j,k}'=0$ \ for \ $1\leq j,k\leq3$
         \ with \ $(j,k)\ne(3,3)$,
  \item $b_{1,1}'=0$ \ and \ $b_{1,1}''=0$,
 \end{enumerate}
 where \ $\delta':=\sqrt{b_{1,1}'b_{2,2}'-(b_{1,2}')^2}$ \ and
 \ $\delta'':=\sqrt{b_{1,1}''b_{2,2}''-(b_{1,2}'')^2}$.
\ In cases \ $(\widetilde{C}1)$, $(\widetilde{C}3)$, $(\widetilde{C}5)$ \ the
 parameters of the convolution \ $\mu'*\mu''$ \ are \ $(a'+a'',B'+B'')$, \ but
 in the other cases it does not necessarily hold (compare with Lemma
 \ref{GAUSSPAR2}).
\end{Thm}

\noindent
\textbf{Proof.} \
First we show necessity, i.e., if \ $\mu'*\mu''$ \ is a Gaussian measure then
 one of the conditions \ ${\rm(\tC1)-(\tC7)}$ \ holds.
Let us denote the parameters of the convolution \ $\mu'*\mu''$ \ by \ $(a,B)$
 \ and we write \ $d_{j,k}:=d_{j,k}^{\pm\lambda}(a,B)$,
 \ $d'_{j,k}:=d_{j,k}^{\pm\lambda}(a',B')$ \ and
 \ $d''_{j,k}:=d_{j,k}^{\pm\lambda}(a'',B'')$ \ for \ $1\leq j,k\leq3$ \ as in
 Theorem \ref{CONV}.
If \ $b_{1,1}'>0$ \ and \ $b_{1,1}''>0$, \ we can easily prove that
 $$
  \frac{b_{1,2}}{b_{1,1}}=\frac{b'_{1,2}}{b'_{1,1}}
   =\frac{b''_{1,2}}{b''_{1,1}},\qquad
  \frac{b_{2,2}}{b_{1,1}}=\frac{b'_{2,2}}{b'_{1,1}}
   =\frac{b''_{2,2}}{b''_{1,1}},
 $$
and \ $d_{2,2}'+d_{1,1}''\in\RR$ \ as in Pap \cite[Theorem 7.3]{PAP02}.
This implies that there exists \ $\varrho>0$ \ such that
 \ $b_{j,k}''=\varrho b_{j,k}'$ \ for \ $1\leq j,k\leq2$, \ i.e.,
 \ ${\rm (\tC2)}$ \ holds.

When \ $b_{1,1}'>0$, \ $\delta'>0$ \ and \ $b_{1,1}''>0$, \ $\delta''>0$, \ we
 show that \ ${\rm (\widetilde{C}1)}$ \ holds.
To derive this it is sufficient to show that \ $b_{1,3}''=\varrho b_{1,3}'$,
 \ $b_{2,3}''=\varrho b_{2,3}'$, \ $a_1''=\varrho a_1'$ \ and
 \ $a_2''=\varrho a_2'$.
\ Using Theorem \ref{CONV} we obtain
\renewcommand{\labelenumi}{\textup{(\roman{enumi})}}
 \begin{enumerate}
  \item $(d_{2,2}'+d_{1,1}'')(\Re d_{1,3}'-\Re d_{1,3})
         =d_{1,2}'(\Re d_{1,3}''+\Re d_{2,3}')$,
  \item $(d_{2,2}'+d_{1,1}'')(\Re d_{2,3}''-\Re d_{2,3})
         =d_{1,2}''(\Re d_{1,3}''+\Re d_{2,3}')$,
  \item $(d_{2,2}'+d_{1,1}'')(\Im d_{1,3}'-\Im d_{1,3})
         =d_{1,2}'(\Im d_{1,3}''+\Im d_{2,3}')$,
  \item $(d_{2,2}'+d_{1,1}'')(\Im d_{2,3}'-\Im d_{2,3})
         =d_{1,2}''(\Im d_{1,3}''+\Im d_{2,3}')$.
 \end{enumerate}
Let us denote
 \ $\delta_1':=b_{1,1}'b_{2,3}'-b_{1,2}'b_{1,3}'$,
 \ $\delta_1'':=b_{1,1}''b_{2,3}''-b_{1,2}''b_{1,3}''$,
 \ $\delta_2':=a_1'b_{1,2}'-a_2'b_{1,1}'$,
 \ $\delta_2'':=a_1''b_{1,2}''-a_2''b_{1,1}''$.
\ Summing up \ ${\rm (iii)}$ \ and \ ${\rm (iv)}$ \ we have
 $$
  (d_{2,2}'+d_{1,1}'')(\Im d_{1,3}'+\Im d_{2,3}''-\Im d_{1,3}-\Im d_{2,3})
  =(d_{1,2}'+d_{1,2}'')(\Im d_{1,3}''+\Im d_{2,3}').
 $$
Using the definition of \ $d_{j,k},\,d_{j,k}',\,d_{j,k}''$
 \ $(1\leq j,k\leq 3)$ \ we get
 \begin{align*}
  &\big(\coth(\lambda\delta')+\coth(\lambda\delta'')\big)
   \bigg(\frac{b_{1,3}'}{b_{1,1}'}
         +\frac{\delta_2'}{\lambda b_{1,1}'\delta'\coth(\lambda\delta'/2)}
         -\frac{b_{1,3}''}{b_{1,1}''}
        +\frac{\delta_2''}{\lambda b_{1,1}''\delta''\coth(\lambda\delta''/2)}\\
  &\phantom{\big(\coth(\lambda\delta')+\coth(\lambda\delta'')\big)\bigg(}
   -\frac{2\delta_2}{\lambda b_{1,1}\delta\coth(\lambda\delta/2)}\bigg)\\
  &=-\bigg(\frac{1}{\sinh(\lambda\delta')}
           +\frac{1}{\sinh(\lambda\delta'')}\bigg)
     \bigg(\frac{b_{1,3}''}{b_{1,1}''}
           +\frac{\delta_2''}
                 {\lambda b_{1,1}''\delta''\coth(\lambda\delta''/2)}
           -\frac{b_{1,3}'}{b_{1,1}'}
           +\frac{\delta_2'}
                 {\lambda b_{1,1}'\delta'\coth(\lambda\delta'/2)}\bigg).
 \end{align*}
An easy calculation shows that
 \begin{align*}
  &\bigg(\frac{b_{1,3}'}{b_{1,1}'}-\frac{b_{1,3}''}{b_{1,1}''}\bigg)
   \lambda\sinh(\lambda\delta'/2)\sinh(\lambda\delta''/2)\\
  &=\bigg(\frac{1}{\delta'+\delta''}\Big(a_1\frac{b_{1,2}}{b_{1,1}}-a_2\Big)
          -\frac{1}{\delta'}\Big(a_1'\frac{b_{1,2}'}{b_{1,1}'}-a_2'\Big)\bigg)
    \sinh(\lambda\delta'/2)\cosh(\lambda\delta''/2)\\
  &\phantom{=\:}
   +\bigg(\frac{1}{\delta'+\delta''}\Big(a_1\frac{b_{1,2}}{b_{1,1}}-a_2\Big)
          -\frac{1}{\delta''}
           \Big(a_1''\frac{b_{1,2}''}{b_{1,1}''}-a_2''\Big)\bigg)
    \cosh(\lambda\delta'/2)\sinh(\lambda\delta''/2)
 \end{align*}
 for all \ $\lambda>0$.
\ We show that the functions
 \ $\lambda\sinh(\lambda\delta'/2)\sinh(\lambda\delta''/2)$,
 \ $\sinh(\lambda\delta'/2)\cosh(\lambda\delta''/2)$ \ and
 \ $\cosh(\lambda\delta'/2)\sinh(\lambda\delta''/2)$ \ $(\lambda>0)$ \ are
 linearly independent.
We have
 \begin{align*}
  \lambda\sinh(\lambda\delta'/2)\sinh(\lambda\delta''/2)
  &=\lambda
    \big(e^{\lambda(\delta'+\delta'')/2}-e^{\lambda(\delta''-\delta')/2}
         -e^{\lambda(\delta'-\delta'')/2}
         +e^{-\lambda(\delta'+\delta'')/2}\big)/4,\\
  \sinh(\lambda\delta'/2)\cosh(\lambda\delta''/2)
  &=\big(e^{\lambda(\delta'+\delta'')/2}+e^{\lambda(\delta'-\delta'')/2}
         -e^{\lambda(\delta''-\delta')/2}
         -e^{-\lambda(\delta'+\delta'')/2}\big)/4,\\
  \cosh(\lambda\delta'/2)\sinh(\lambda\delta''/2)
  =&\big(e^{\lambda(\delta'+\delta'')/2}-e^{\lambda(\delta'-\delta'')/2}
         +e^{\lambda(\delta''-\delta')/2}
         -e^{-\lambda(\delta'+\delta'')/2}\big)/4.
 \end{align*}
The linear independence of these functions follows from the following fact: if
 \ $c_1,\dots,c_n$ \ are pairwise different complex numbers \ and
 \ $Q_1,\dots,Q_n$ \ are complex polynomials such that
 \ $\sum_{j=1}^{n}Q_j(\lambda)e^{c_j\lambda}=0$ \ for all \ $\lambda>0$ \ then
 \ $Q_1=\cdots =Q_n=0$.
\ Hence we get
 \begin{equation}
  \label{SCHREPPRELTOOL2}
  \frac{b_{1,3}'}{b_{1,1}'}-\frac{b_{1,3}''}{b_{1,1}''}=0,\qquad
  \frac{1}{\delta'+\delta''}\big(a_1\frac{b_{1,2}}{b_{1,1}}-a_2\big)
  =\frac{1}{\delta'}\big(a_1'\frac{b_{1,2}'}{b_{1,1}'}-a_2'\big)
  =\frac{1}{\delta''}\big(a_1''\frac{b_{1,2}''}{b_{1,1}''}-a_2''\big).
 \end{equation}
Subtracting the equation \ ${\rm(i)}$ \ from \ ${\rm(ii)}$ \ we get
 $$
  (d_{2,2}'+d_{1,1}'')(\Re d_{1,3}'-\Re d_{2,3}''-\Re d_{1,3}+\Re d_{2,3})
  =(d_{1,2}'-d_{1,2}'')(\Re d_{1,3}''+\Re d_{2,3}').
 $$
Using again the definition of \ $d_{j,k}$, $d_{j,k}'$, $d_{j,k}''$
 \ $(1\leq j,k\leq 3)$ \ we obtain
 \begin{align*}
  &\big(\coth(\lambda\delta')+\coth(\lambda\delta'')\big)
   \bigg(\frac{a_1'}{\sqrt{\lambda}b_{1,1}'}
         +\frac{a_1''}{\sqrt{\lambda}b_{1,1}''}
         -\frac{2a_1}{\sqrt{\lambda}b_{1,1}}\\
  &\phantom{\big(\coth(\lambda\delta')+\coth(\lambda\delta'')\big)\bigg(}
         +\frac{\sqrt{\lambda}\delta_1'}
               {b_{1,1}'\delta'\coth(\lambda\delta'/2)}
         -\frac{\sqrt{\lambda}\delta_1''}
               {b_{1,1}''\delta''\coth(\lambda\delta''/2)}\bigg)\\
  &=\bigg(\frac{1}{\sinh(\lambda\delta'')}
          -\frac{1}{\sinh(\lambda\delta')}\bigg)
    \bigg(\frac{a_1''}{\sqrt{\lambda}b_{1,1}''}
          -\frac{a_1'}{\sqrt{\lambda}b_{1,1}'}
          +\frac{\sqrt{\lambda}\delta_1'}
                {b_{1,1}'\delta'\coth(\lambda\delta'/2)}
          +\frac{\sqrt{\lambda}\delta_1''}
                {b_{1,1}''\delta''\coth(\lambda\delta''/2)}\bigg).
 \end{align*}
A simple calculation shows that
 \begin{align*}
  &\lambda\big(1+\tanh(\lambda\delta'/2)\tanh(\lambda\delta''/2)\big)
   \bigg(\frac{\delta_1'}{\delta'b_{1,1}'}
         -\frac{\delta_1''}{\delta''b_{1,1}''}\bigg)\\
  &=\big(\coth(\lambda\delta')+\coth(\lambda\delta'')\big)
    \bigg(2\frac{a_1}{b_{1,1}}-\frac{a_1'}{b_{1,1}'}
          -\frac{a_1''}{b_{1,1}''}\bigg)
    +\bigg(\frac{1}{\sinh(\lambda\delta')}
           -\frac{1}{\sinh(\lambda\delta'')}\bigg)
     \bigg(\frac{a_1'}{b_{1,1}'}-\frac{a_1''}{b_{1,1}''}\bigg).
 \end{align*}
It can be easily checked that the functions
 \ $\lambda\big(1+\tanh(\lambda\delta'/2)\tanh(\lambda\delta''/2)\big)$,
 \ $\coth(\lambda\delta')+\coth(\lambda\delta'')$ \ and
 \ $(\sinh(\lambda\delta'))^{-1}-(\sinh(\lambda\delta''))^{-1}$ \ $(\lambda>0)$
 \ are linearly independent.
Hence we have
 \begin{equation}
  \label{SCHREPPRELTOOL4}
  \frac{a_1'}{b_{1,1}'}-\frac{a_1''}{b_{1,1}''}=0,\qquad
  2\frac{a_1}{b_{1,1}}-\frac{a_1'}{b_{1,1}'}-\frac{a_1''}{b_{1,1}''}=0,\qquad
  \frac{\delta_1'}{\delta'b_{1,1}'}=\frac{\delta_1''}{\delta''b_{1,1}''}.
 \end{equation}
Taking into account \eqref{SCHREPPRELTOOL2} and \eqref{SCHREPPRELTOOL4}, we
 conclude that \ ${\rm(\tC1)}$ \ holds.
Using Lemma \ref{GAUSSPAR2} it turns out that in this case \ $a=a'+a''$ \ and
 \ $B=B'+B''$.

If \ $b_{1,1}'>0$, \ $\delta'>0$ \ and \ $b_{1,1}''>0$, \
$\delta''=0$ \ we show that \ $\mu'*\mu''$ \ can not be a Gaussian
measure. Our proof goes along the lines of the proof Theorem 7.3 in
Pap \cite{PAP02}. Since the proof given in Pap \cite{PAP02} contains
a mistake we write down the details. Suppose that, on the contrary,
\ $\mu'*\mu''$ \ is a Gaussian measure on \ $\HH$ \ with parameters
\ $(a,B).$ \ By Lemma \ref{GAUSSPAR2}, we have \
$b_{1,1}=b_{1,1}'+b_{1,1}''$, \ hence \ $b_{1,1}>0.$ \ By Theorem
\ref{SREP}, we have
 \ $\left(\mu'*\mu''\right)\ssphat(\pi_{\pm\lambda})$ \ is an integral operator.
 Using Theorem \ref{CONV} we obtain
 \begin{align}\label{SEGED_NO_GAUSS1}
   &d_{1,1}=d_{1,1}'-\frac{(d_{1,2}')^2}{d_{2,2}'+d_{1,1}''},\\\label{SEGED_NO_GAUSS2}
   &d_{2,2}=d_{2,2}''-\frac{(d_{1,2}'')^2}{d_{2,2}'+d_{1,1}''}.
 \end{align}
 We show that \ $d_{2,2}'+d_{1,1}''\in\RR$ \ and \ $\frac{b_{1,2}'}{b_{1,1}'}
  =\frac{b_{1,2}''}{b_{1,1}''}.$ \
 (The derivations of these two facts are not correct in the proof
  of Theorem 7.3 in Pap \cite{PAP02}.)
By Theorem \ref{SREP}, we have
$$
  \Im(d_{2,2}'+d_{1,1}'')
   =\mp\left(\frac{b_{1,2}'}{b_{1,1}'}-\frac{b_{1,2}''}{b_{1,1}''}\right)
   =-\Im(d_{1,1}'+d_{2,2}'').
$$
Using that \ $\Im(d_{1,1}+d_{2,2})=0,$ \ by \eqref{SEGED_NO_GAUSS1}
and \eqref{SEGED_NO_GAUSS2} we get
\begin{align*}
 0&=\pm\left(\frac{b_{1,2}'}{b_{1,1}'}-\frac{b_{1,2}''}{b_{1,1}''}\right)
   -\Im\left(\frac{(d_{1,2}')^2+(d_{1,2}'')^2}{d_{2,2}'+d_{1,1}''}\right)\\
  &=\pm\left(\frac{b_{1,2}'}{b_{1,1}'}-\frac{b_{1,2}''}{b_{1,1}''}\right)
   \mp\frac{(d_{1,2}')^2+(d_{1,2}'')^2}{\vert d_{2,2}'+d_{1,1}''\vert^2}
      \left(\frac{b_{1,2}'}{b_{1,1}'}-\frac{b_{1,2}''}{b_{1,1}''}\right).
\end{align*}
 Hence
 $$
   \Big(\vert d_{2,2}'+d_{1,1}''\vert^2-(d_{1,2}')^2-(d_{1,2}'')^2\Big)
       \left(\frac{b_{1,2}'}{b_{1,1}'}-\frac{b_{1,2}''}{b_{1,1}''}\right)=0.
 $$
 Then
 \begin{align*}
  \vert d_{2,2}'+d_{1,1}''\vert^2-(d_{1,2}')^2-(d_{1,2}'')^2
    &=\left\vert\frac{\delta'\coth(\lambda\delta')\mp ib_{1,2}'}{b_{1,1}'}
           +\frac{\lambda^{-1}\pm ib_{1,2}''}{b_{1,1}''}\right\vert^2
     -\frac{(\delta')^2}{(b_{1,1}')^2\sinh^2(\lambda\delta')}\\
    &\phantom{=}
       -\frac{1}{\lambda^2(b_{1,1}'')^2}\\
    &=\frac{(\delta')^2}{(b_{1,1}')^2}
       +\frac{2\delta'\coth(\lambda\delta')}{\lambda b_{1,1}'b_{1,1}''}
       +\left(\frac{b_{1,2}'}{b_{1,1}'}-\frac{b_{1,2}''}{b_{1,1}''}\right)^2
     >0.
 \end{align*}
It yields that \
 $\frac{b_{1,2}'}{b_{1,1}'}=\frac{b_{1,2}''}{b_{1,1}''}.$ \
 Particularly, \ $d_{2,2}'+d_{1,1}''\in\RR.$ \ Rewrite \eqref{SEGED_NO_GAUSS1}
and \eqref{SEGED_NO_GAUSS2} in the form
\begin{align*}
  &(d_{1,1}'-d_{1,1})(d_{2,2}'+d_{1,1}'')=(d_{1,2}')^2,\\
  &(d_{2,2}''-d_{2,2})(d_{2,2}'+d_{1,1}'')=(d_{1,2}'')^2.
\end{align*}
 It follows that
 $$
   (d_{1,1}'-d_{2,2}''-d_{1,1}+d_{2,2})(d_{2,2}'+d_{1,1}'')
      =(d_{1,2}')^2-(d_{1,2}'')^2.
 $$
 Using that \ $d_{2,2}'+d_{1,1}''\in\RR$ \ and \ $\Re(d_{1,1}-d_{2,2})=0,$ \
 taking real parts we get
 $$
  (\Re(d_{1,1}')-\Re(d_{2,2}''))(d_{2,2}'+d_{1,1}'')
      =(d_{1,2}')^2-(d_{1,2}'')^2.
 $$
 Thus
\begin{align*}
  \left(\frac{\delta'\coth(\lambda\delta')}{b_{1,1}'}
         -\frac{1}{\lambda b_{1,1}''}\right)
  \left(\frac{\delta'\coth(\lambda\delta')}{b_{1,1}'}
         +\frac{1}{\lambda b_{1,1}''}\right)
   =\frac{(\delta')^2}{(b_{1,1}')^2\sinh^2(\lambda\delta')}
    -\frac{1}{\lambda^2 (b_{1,1}'')^2}.
\end{align*}
 From this we conclude
 $$
    \frac{(\delta')^2\coth^2(\lambda\delta')}{(b_{1,1}')^2}
      -\frac{1}{\lambda^2 (b_{1,1}'')^2}
    =\frac{(\delta')^2}{(b_{1,1}')^2\sinh^2(\lambda\delta')}
      -\frac{1}{\lambda^2 (b_{1,1}'')^2}
     \quad\Longleftrightarrow\quad
     \cosh(\lambda\delta')=1,
 $$
thus \ $\delta'=0,$ \ which leads to a contradiction.

If \ $b_{1,1}'>0$, \ $\delta'>0$, \ and \ $b_{1,1}''=0$ \ we show that
 \ ${\rm(\tC3)}$ \ holds.
The symmetry and positive semi--definiteness of the matrix \ $B''$ \ imply
 \ $b_{1,2}''=b_{1,3}''=0$.
\ Lemma \ref{GAUSSPAR2} yields that \ $b_{1,1}=b_{1,1}'+b_{1,1}''>0$.
\ Hence Theorem \ref{SREP} implies that
 \ $\left(\mu'*\mu''\right)\ssphat(\pi_{\pm\lambda})$ \ is an integral operator
 and \ $\Im(d_{1,1}+d_{2,2})=0$ \ holds.
By Theorem \ref{SREP} and Theorem \ref{CONV} we obtain
 \ $\Im(d_{1,1}+d_{2,2})
    =\Im(d_{1,1}'+d_{2,2}'+\lambda b_{2,2}'')=\Im(\lambda b_{2,2}'')$.
\ Thus \ $b_{2,2}''=0$, \ which implies that \ $b_{2,3}''=0$ \ and
 \ $\delta=\delta'>0$.
\ Using again Theorem \ref{CONV} we get
 \begin{align}
  \label{SCHREPPRELTOOL5}
  d_{1,3}&=d_{1,3}'-\sqrt{\lambda}a_1''d_{1,2}'\\\label{SCHREPPRELTOOL6}
  d_{2,3}&=d_{2,3}'-\sqrt{\lambda}a_1''d_{2,2}'\mp i\sqrt{\lambda}a_2''.
 \end{align}
Taking the real part of the difference of equations \eqref{SCHREPPRELTOOL5} and
 \eqref{SCHREPPRELTOOL6} we have
 \begin{equation}
  \label{SCHREPPRELTOOL7}
  2\bigg(\frac{a_1}{b_{1,1}}-\frac{a_1'}{b_{1,1}'}\bigg)
  =\lambda\delta'\frac{a_1''}{b_{1,1}'}
   \bigg(\frac{1+\cosh(\lambda\delta')}{\sinh(\lambda\delta')}\bigg).
 \end{equation}
Since \eqref{SCHREPPRELTOOL7} is valid for all \ $\lambda>0$, \ we have
 \ $a_1''=0$.
\ Taking the imaginary part of \eqref{SCHREPPRELTOOL6} and using the fact that
 \ $a_1''=0$ \ we get
 \begin{equation}
  \label{SCHREPPRELTOOL8}
  a_2''\bigg(1-\frac{1}{\lambda\delta'\coth(\lambda\delta'/2)}\bigg)
  =\frac{b_{1,3}}{b_{1,1}}-\frac{b_{1,3}'}{b_{1,1}'}=0.
 \end{equation}
Since \eqref{SCHREPPRELTOOL8} is valid for all \ $\lambda>0$, \ we get
 \ $a_2''=0$, \ so \ ${\rm (\widetilde{C}3)}$ \ holds.
If \ $b_{1,1}'>0$, \ $\delta'=0$ \ and \ $b_{1,1}''=0$ \ a similar argument
 shows that \ ${\rm (\widetilde{C}4)}$ \ holds.

The aim of the following discussion is to show the converse.
Suppose that \ ${\rm (\widetilde{C}1)}$ \ holds.
We prove that the convolution \ $\mu'*\mu''$ \ is a Gaussian measure on
 \ $\HH$ \ with parameters \ $(a'+a'',B'+B'')$.
\ By Theorem \ref{CONV}, \ the Fourier transform
 \ $\left(\mu'*\mu''\right)\ssphat(\chi_{\alpha,\beta})$ \ equals the Fourier
 transform of a Gaussian measure with parameters \ $(a'+a'',B'+B'')$ \ at the
 representation \ $\chi_{\alpha,\beta}$ \ for all \ $\alpha,\beta>0$.
\ Since \ $b_{1,1}'+b_{1,1}''>0$, \ the Fourier transform
 \ $\left(\mu'*\mu''\right)\ssphat(\pi_{\pm\lambda})$ \ is an integral operator
 on \ $L^2(\RR)$ \ with kernel function \ $K_{\pm\lambda}$ \ given in Theorem
 \ref{CONV} for all \ $\lambda>0$.
\ All we have to show is that \ $C=C_{\pm\lambda}(B'+B'')$ \ and \
 $V=D_{\pm\lambda}(a'+a'',B'+B'')
  =(d_{j,k}^{\pm\lambda}(a'+a'',B'+B''))_{1\sleq j,k\sleq3}$.
\ We have
\begin{align*}
  d_{2,2}'+d_{1,1}''
  =\frac{\delta'\sinh\big(\lambda(1+\varrho)\delta'\big)}
        {b_{1,1}'\sinh(\lambda\delta')\sinh(\lambda\varrho\delta')},
 \end{align*}
 hence using Theorem \ref{CONV} we obtain
 $$
  C=\sqrt{\frac{\delta'}{2\pi b_{1,1}'\sinh(\lambda(1+\varrho)\delta')}}
   =C_{\pm\lambda}(B'+B'').
 $$
Let \ $(\mu_t)_{t\sgeq 0}$ \ be a Gaussian semigroup such that \ $\mu_1$ \ is a
 Gaussian measure with parameters \ $(a',B')$.
\ By the help of the semigroup property  we have
 \ $\mu_1*\mu_\varrho=\mu_{1+\varrho}$.
\ Taking into account that \ $a_3'$ \ and \ $b_{3,3}'$ \ appear only in
 \ $d_{3,3}^{\pm\lambda}(a',B')$ \ (see Theorem \ref{SREP}) and the fact that
 \ $\mu_t$ \ is a Gaussian measure with parameters \ $(ta',tB')$ \ for all
 \ $t\geq 0,$ \ Theorem \ref{SREP} and Theorem \ref{CONV} give us
 $$
  v_{j,k}=d_{j,k}^{\pm\lambda}(a'+a'',B'+B'').
 $$
 for \ $1\leq j,k\leq 3$ \ with \ $(j,k)\ne(3,3)$.
\ So we have to check only that
 \ $v_{3,3}=d_{3,3}^{\pm\lambda}(a'+a'',B'+B'')$.
\ By the help of Theorem \ref{CONV} we get
 \begin{equation}\label{SCHREPPRELTOOL9}
  v_{3,3}
  =d_{3,3}'+d_{3,3}''-\frac{1}{d_{2,2}'+d_{1,1}''}(d_{3,2}'+d_{3,1}'')^2.
 \end{equation}
Calculating the real and imaginary part of \eqref{SCHREPPRELTOOL9} one can
 easily check that \ $v_{3,3}=d_{3,3}^{\pm\lambda}(a'+a'',B'+B'')$ \ is valid.

Now suppose that \ ${\rm (\widetilde{C}2)}$ \ holds.
Using the parameters of \ $\mu'$ \ and \ $\mu''$, \ define a vector
 \ $a=(a_i)_{1\sleq i\sleq 3}$ \ and a matrix
 \ $B=(b_{i,j})_{1\sleq i,j\sleq 3}$, \ as in Lemma \ref{GAUSSPAR2}.
We show that the convolution \ $\mu:=\mu'*\mu''$ \ is a Gaussian measure on
 \ $\HH$ \ with parameters \ $(a,B)$.
\ An easy calculation shows that the Fourier transforms of \ $\mu'*\mu''$ \ and
 \ $\mu$ \ at the one--dimensional representations coincide.
Concerning the Fourier transforms at the Schr\"odinger representations, as in
 case of \ ${\rm (\tC1),}$ \ all we have to prove is that
 $$
  C_{\pm\lambda}(B)
  =C_{\pm\lambda}(B')C_{\pm\lambda}(B'')\sqrt{\frac{2\pi}{d_{2,2}'+d_{1,1}''}}
 $$
 and \ $V=D_{\pm\lambda}(a'+a'',B'+B'')$.
\ Using Theorem \ref{SREP} we have
 \begin{align*}
  \frac{1}{\sqrt{2\pi\lambda b_{1,1}'}}
  \frac{1}{\sqrt{2\pi\lambda b_{1,1}''}}
  \sqrt{\frac{2\pi}{\frac{1}{\lambda b_{1,1}'}+\frac{1}{\lambda b_{1,1}''}
                    \pm i\big(\frac{b_{1,2}''}{b_{1,1}''}
                          -\frac{b_{1,2}'}{b_{1,1}'}\big)}}
  =\frac{1}{\sqrt{2\pi\lambda(b_{1,1}'+b_{1,1}'')}}
  =\frac{1}{\sqrt{2\pi\lambda b_{1,1}}},
 \end{align*}
 since \ $b_{1,2}''/b_{1,1}''=b_{1,2}'/b_{1,1}'=\varrho$.
\ Using similar arguments one can also easily check that
 \ $V=D_{\pm\lambda}(a'+a'',B'+B'')$ \ holds.
We note that in this case the parameters of \ $\mu'*\mu''$ \ is not the sum of
 the parameters of \ $\mu'$ \ and \ $\mu''$.

Suppose that \ ${\rm (\tC3)}$ \ holds.
Proposition \ref{CENTRALGAUSS} gives us that the convolution \ $\mu'*\mu''$
 \ is a Gaussian measure on \ $\HH$ \ with parameters \ $(a'+a'',B'+B'')$.
\ In cases \ ${\rm(\tC4)}$, ${\rm(\tC5)}$, ${\rm(\tC6)}$, ${\rm(\tC7)}$ \ we
 can argue as in cases \ ${\rm (\widetilde{C}2),\,(\widetilde{C}3).}$
\ Consequently, the proof is complete.
\proofend

For the proof of Theorem \ref{C} we need the following lemma about the support
 of a Gaussian measure on \ $\HH$.

\begin{Lem}\label{SUPPEU} \
Let \ $\mu$ \ be a Gaussian measure on \ $\HH$ \ with parameters \ $(a,B)$
 \ such that \ $b_{1,1}b_{2,2}-b_{1,2}^2=0$.
\ Let \ $Y_0\in\cH$ \ be defined as in Section \ref{PREL}.
If \ $\rank(B)=2$ \ then
 \ $\supp(\mu)=\exp\big(Y_0+\RR\cdot U+\RR\cdot X_3\big)$, \ where
 $$
  U:=\begin{cases}
      b_{1,1}X_1+b_{2,1}X_2 & \text{if \ $b_{1,1}>0$,}\\
      b_{2,2}X_2 & \text{if \ $b_{1,1}=0$ \ and \ $b_{2,2}>0$.}
     \end{cases}
 $$
If \ $\rank(B)=1$ \ then
 \ $\supp(\mu)=\exp\big(Y_0+\RR\cdot U+\RR\cdot [Y_0,U]\big)$, \ where
 $$U:=\begin{cases}
       b_{1,1}X_1+b_{2,1}X_2+b_{3,1}X_3 & \text{if \ $b_{1,1}>0$,}\\
       b_{2,2}X_2+b_{3,2}X_3 & \text{if \ $b_{1,1}=0$ \ and \ $b_{2,2}>0$,}\\
       b_{3,3}X_3 & \text{if \ $b_{1,1}=b_{2,2}=0$ \ and \ $b_{3,3}>0$.}
     \end{cases}$$
If \ $\rank(B)=0$ \ then \ $\supp(\mu)=\exp(Y_0)$.
\end{Lem}

\noindent
\textbf{Proof.} \
We apply \ ${\rm (iii)-(v)}$ \ of Lemma \ref{GAUSSSUPPORT}, respectively.
If \ $\rank(B)=2$ \ then one can check that \ $\cL(Y_1,Y_2)=\cL(U,X_3)$.
\ If \ $\rank(B)=1$ \ then \ $\cL(Y_1)=\cL(U)$.
\proofend

\noindent
\textbf{Proof of Theorem \ref{C}.} \
First we prove that if one of the conditions \ ${\rm (C1)}$ \ and
 \ ${\rm (C2)}$ \ holds then one of the conditions \ ${\rm (\tC1)-(\tC7)}$ \ in
 Theorem \ref{SCHREPPREL} is valid, which implies that the convolution
 \ $\mu'*\mu''$ \ is a Gaussian measure on \ $\HH$.

Suppose that \ ${\rm(C1)}$ \ holds.
Lemma \ref{GAUSSSUPPORT} implies \ $\delta'=\delta''=0$.

If \ $b_{1,1}'=b_{1,1}''=0$ \ then \ ${\rm (\widetilde{C}7)}$ \ holds.

If \ $b_{1,1}'>0$, \ $\delta'=0$ \ and \ $b_{1,1}''=0$, \ $\delta''=0$ \ we
 show that \ ${\rm(\tC4)}$ \ holds.
It is sufficient to show that \ $b_{2,2}''=0$.
\ Suppose that, on the contrary, \ $b_{2,2}''\ne 0$.
\ When \ $\rank(B')=\rank(B'')=2$, \ by the help of Lemma \ref{SUPPEU}, we get
 $$
  \supp(\mu')=\exp\big(Y_0'+\RR\cdot U'+\RR\cdot X_3\big),\qquad
  \supp(\mu'')=\exp\big(Y_0''+\RR\cdot U''+\RR\cdot X_3\big),
 $$
 where \ $U'=b_{1,1}'X_1+b_{2,1}'X_2$ \ and \ $U''=b_{2,2}''X_2$.
\ Since in this case \ $\supp(\mu')$ \ and \ $\supp(\mu'')$ \ are contained in
 ``Euclidean cosets'' of the same 2--dimensional Abelian subgroup of \ $\HH$,
 \ we obtain that \ $\cL(U',X_3)=\cL(U'',X_3)$.
\ From this we conclude \ $b_{1,1}'=0$, \ which leads to a contradiction.
When \ $\rank(B')=1$, \ $\rank(B'')=2$ \ and in other cases one can argue
 similarly, so \ ${\rm (\tC4)}$ \ holds.

If \ $b_{1,1}'=0$, \ $\delta'=0$ \ and \ $b_{1,1}''>0$, \ $\delta''=0$ \ the
 same argument shows that \ ${\rm(\tC6)}$ \ holds.

If \ $b_{1,1}'>0$, \ $\delta'=0$ \ and \ $b_{1,1}''>0$, \ $\delta''=0$ \ we
 show that \ ${\rm(\tC2)}$ \ holds.
When \ $\rank(B')=\rank(B'')=2$, \ Lemma \ref{SUPPEU} implies that
 $$
  \supp(\mu')=\exp\big(Y_0'+\RR\cdot U'+\RR\cdot X_3\big),\qquad
  \supp(\mu'')=\exp\big(Y_0''+\RR\cdot U''+\RR\cdot X_3\big),
 $$
 where \ $U'=b_{1,1}'X_1+b_{2,1}'X_2$ \ and \ $U''=b_{1,1}''X_1+b_{2,1}''X_2$.
\ Condition \ ${\rm(C1)}$ \ yields that \ $\cL(U',X_3)=\cL(U'',X_3)$, \ hence
 we have \ $b_{2,1}''b_{1,1}'=b_{2,1}'b_{1,1}''$.
\ Since \ $\delta'=\delta''=0$ \ we get
 \ $b_{2,2}''b_{1,1}'=b_{2,2}'b_{1,1}''$.
\ Thus \ ${\rm (\widetilde{C}2)}$ \ holds with
 \ $\varrho:=b_{1,1}''/b_{1,1}'$.
\ When \ $\rank(B')=\rank(B'')=1$, \ Lemma \ref{SUPPEU} implies that
 $$
  \supp(\mu')=\exp\big(Y_0'+\RR\cdot U'+\RR\cdot [Y_0',U']\big),\qquad
  \supp(\mu'')=\exp\big(Y_0''+\RR\cdot U''+\RR\cdot [Y_0'',U'']\big),
 $$
 where \ $U'=b_{1,1}'X_1+b_{2,1}'X_2+b_{3,1}'X_3$ \ and
 \ $U''=b_{1,1}''X_1+b_{2,1}''X_2+b_{3,1}''X_3$.
\ Condition \ ${\rm(C1)}$ \ yields \ $\cL(U',[Y_0',U'])=\cL(U'',[Y_0'',U''])$,
 \ hence \ $\cL(b_{1,1}'X_1+b_{2,1}'X_2)=\cL(b_{1,1}''X_1+b_{2,1}''X_2)$.
\ It can be easily checked that \ ${\rm (\widetilde{C}2)}$ \ holds with
 \ $\varrho:=b_{1,1}''/b_{1,1}'$.
\ When \ $\rank(B')=1$, \ $\rank(B'')=2$ \ or \ $\rank(B')=2$, \ $\rank(B'')=1$
 \ we also have \ ${\rm(\tC2)}$ \ holds.

Suppose that \ ${\rm(C2)}$ \ holds (i.e., $\mu'=\mu_{t'}$,
 \ $\mu''=\mu_{t''}*\nu$ \ or \ $\mu'=\mu_{t'}*\nu$, \ $\mu''=\mu_{t''}$ \ with
 appropriate nonnegative real numbers \ $t'$, $t''$ \ and a Gaussian measure
 \ $\nu$ \ with support contained in the center of \ $\HH$).
Then we have
 \begin{align*}
  \mu'*\mu''=\mu_{t'}*\mu_{t''}*\nu=\mu_{t'+t''}*\nu\qquad{\rm or}\qquad
  \mu'*\mu''=\mu_{t'}*\nu*\mu_{t''}=\mu_{t'+t''}*\nu.
 \end{align*}
Remark \ref{CENTRALGAUSSREM} and Proposition \ref{CENTRALGAUSS}
yield that
 \ $\mu'*\mu''$ \ is a Gaussian measure on \ $\HH$.

Conversely, suppose that \ $\mu'*\mu''$ \ is a Gaussian measure on \ $\HH$.
\ Then by Theorem \ref{SCHREPPREL}, one of the conditions
 \ ${\rm(\tC1)-(\tC7)}$ \ holds.
We show that then one of the conditions \ ${\rm(C1)}$ \ and \ ${\rm(C2)}$ \ is
 valid.

Suppose that \ ${\rm(\tC1)}$ \ holds.
If \ $b_{3,3}''-\varrho b_{3,3}'\geq 0$ \ then let \ $(\alpha_t')_{t\sgeq0}$
 \ be a Gaussian semigroup such that \ $\alpha_1'=\mu'$ \ and let \ $\nu$ \ be
 a Gaussian measure on \ $\HH$ \ with parameters \ $(a_\nu,B_\nu)$ \ such that
 $$
  B_\nu:=\begin{bmatrix}
          0 & 0 & 0 \\
          0 & 0 & 0 \\
          0 & 0 & b_{3,3}''-\varrho b_{3,3}'
        \end{bmatrix},\qquad
  a_\nu:=\begin{bmatrix}
          0 \\
          0  \\
          a_3''-\varrho a_3'
        \end{bmatrix}.
 $$
Remark \ref{CENTRALGAUSSREM} and Proposition \ref{CENTRALGAUSS}
imply that
 \ $\mu''=\alpha_\varrho'*\nu$, \ hence \ ${\rm (C2)}$ \ holds.
If \ $b_{3,3}''-\varrho b_{3,3}'< 0$ \ then let \ $(\alpha_t'')_{t\sgeq 0}$
 \ be a Gaussian semigroup such that \ $\alpha_1''=\mu''$ \ and let \ $\nu$
 \ be a Gaussian measure on \ $\HH$ \ with parameters \ $(a_\nu,B_\nu)$ \ such
 that
 $$
  B_\nu:=\begin{bmatrix}
          0 & 0 & 0 \\
          0 & 0 & 0 \\
          0 & 0 & b_{3,3}'-\varrho^{-1}b_{3,3}''
        \end{bmatrix},\qquad
  a_\nu:=\begin{bmatrix}
          0 \\
          0  \\
          a_3'-\varrho^{-1}a_3'
        \end{bmatrix}.
 $$
Remark \ref{CENTRALGAUSSREM} and Proposition \ref{CENTRALGAUSS}
imply that
 \ $\mu'=\alpha_{1/\varrho}''*\nu$, \ hence \ ${\rm(C2)}$ \ holds.

Suppose that \ ${\rm (\widetilde{C}2)}$ \ holds.
Lemma \ref{SUPPEU} implies that
 $$
  \supp(\mu')\subset\exp\big(Y_0'+\RR\cdot U'+\RR\cdot X_3\big),\qquad
  \supp(\mu'')\subset\exp\big(Y_0''+\RR\cdot U''+\RR\cdot X_3\big),
 $$
 where \ $U'=b_{1,1}'X_1+b_{2,1}'X_2$ \ and \ $U''=b_{1,1}''X_1+b_{2,1}''X_2$.
\ Condition \ ${\rm(\tC2)}$ \ gives us that \ $\cL(U')=\cL(U'')$, \ hence
 \ ${\rm(C1)}$ \ holds.

Suppose that \ ${\rm(\tC3)}$ \ holds.
Let \ $(\alpha_t')_{t\sgeq 0}$ \ be a Gaussian semigroup such that
 \ $\alpha_1'=\mu'$ \ and let \ $\nu$ \ be a Gaussian measure with parameters
 \ $(a_\nu,B_\nu)$ \ such that
 $$
  B_\nu:=\begin{bmatrix}
          0 & 0 & 0 \\
          0 & 0 & 0 \\
          0 & 0 & b_{3,3}''
        \end{bmatrix},\qquad
  a_\nu:=\begin{bmatrix}
          0 \\
          0  \\
          a_3''
        \end{bmatrix}.
 $$
Then we have \ $\mu''=\nu=\alpha_0'*\nu$, \ so \ ${\rm(C2)}$ \ holds.

Suppose that \ ${\rm (\widetilde{C}4)}$ \ holds.
By the help of Lemma \ref{SUPPEU}, we have
 $$
  \supp(\mu')\subset\exp\big(Y_0'+\RR\cdot U'+\RR\cdot X_3\big),\qquad
  \supp(\mu'')\subset\exp\big(Y_0''+\RR\cdot U''\big),
 $$
 where \ $U'=b_{1,1}'X_1+b_{2,1}'X_2$ \ and \ $U''=b_{3,3}''X_3$.
\ Hence the support of \ $\mu'$ \ is contained in
 \ $\exp\big(Y_0'+\RR\cdot U'+\RR\cdot X_3\big)$ \ and the support of \ $\mu''$
 \ is contained in \ $\exp\big(Y_0''+\RR\cdot U'+\RR\cdot X_3\big)$, \ so
 \ ${\rm (C1)}$ \ holds.
Similar arguments show that when \ ${\rm(\tC5)}$ \ holds then \ ${\rm(C2)}$
 \ is valid, and when \ ${\rm(\tC6)}$ \ holds then \ ${\rm(C1)}$ \ is valid.

Suppose that \ ${\rm(\tC7)}$ \ holds.
Using Lemma \ref{SUPPEU}, we have
 $$
  \supp(\mu')\subset\exp\big(Y_0'+\RR\cdot U'+\RR\cdot X_3\big),\qquad
  \supp(\mu'')\subset\exp\big(Y_0''+\RR\cdot U''+\RR\cdot X_3\big),
 $$
 where \ $U'=b_{2,2}'X_2$ \ and \ $U''=b_{2,2}''X_2$, \ so \ ${\rm(C1)}$
 \ holds.
\proofend

\begin{Rem} \
In case of \ ${\rm(C1)}$ \ in Theorem \ref{C}, \ $\mu'$ \ and \ $\mu''$ \ are
 Gaussian measures also in the ``Euclidean sense'' (i.e., considering them as
 measures on \ $\RR^3$), \ but the parameters of the convolution \ $\mu'*\mu''$
 \ are not necessarily the sum of the parameters of \ $\mu'$ \ and \ $\mu''$.
\ In case of \ ${\rm(C2)}$ \ in Theorem \ref{C}, \ $\mu'$ \ and \ $\mu''$ \ are
 not necessarily Gaussian measures in the ''Euclidean sense'', but the
 parameters of the convolution \ $\mu'*\mu''$ \ is the sum of the parameters of
 \ $\mu'$ \ and \ $\mu''$.
\end{Rem}

\begin{Rem}
 We formulate Theorem \ref{C} in the important special case of centered
 Gaussian measures for which the corresponding Gaussian semigroups are stable
 in the sense of Hazod. First we recall that a probability measure \
 $\mu$ \ on \ $\HH$ \ is called centered if
 $$
    \int_{\HH}x_1\;\mu(\dd x)=\int_{\HH}x_2\;\mu(\dd x)=0.
 $$
 A convolution semigroup \ $(\mu_t)_{t\sgeq 0}$ \ on \ $\HH$ \ is called
 centered if \ $\mu_t$ \ is centered for all \ $t\geq0.$ \ For each
 \ $t\geq 0$ \ let \ $d_t$ \ denote the dilation
 $$
    d_t(x)=(tx_1,tx_2,t^2x_3),\quad x\in\HH,\;\;t\geq0.
 $$
 By Hazod \cite[page 229]{HAZ01}, a Gaussian semigroup
 \ $(\mu_t)_{t\sgeq 0}$ \ is centered and stable in the sense that
 \ $\mu_t=d_{\sqrt{t}}\mu_1,$ $t\geq 0$ \ (Hazod stability) if and only if
  its infinitesimal generator has the form
 \begin{align}\label{STABIL}
   a_3X_3+\frac{1}{2}\sum_{i=1}^2\sum_{j=1}^2b_{i,j}X_iX_j.
 \end{align}
 Wehn \cite{WEH62} proved the following central limit theorem. Let \
 $\vert.\vert$ \ be a fixed homogeneous norm on \ $\HH$ \ and  let
 us consider a centered probability measure \ $\mu$ \ on \ $\HH.$ \
 If \ $\int_{\HH}\vert x\vert^2\;\mu(\dd x)<+\infty,$ \ then
 \ $\big(d_{1/\sqrt{n}}(\mu^{*n})\big)_{n\sgeq1}$ \ converges
 towards \ $\nu$ \ weakly, where \ $\nu$ \ is a Gaussian measure
 on \ $\HH$ \ such that the corresponding Gaussian semigroup has
 infinitesimal generator \eqref{STABIL}.

 For centered and stable Gaussian measures Theorem \ref{C} has
 the following form.

{\sl Let \ $\mu'$ \ and \ $\mu''$ \ be Gaussian measures on \ $\HH$
\ such that the corresponding Gaussian semigroups have
 infinitesimal generators
 $$
  a_3'X_3+\frac{1}{2}\sum_{i=1}^2\sum_{j=1}^2b_{i,j}'X_iX_j
   \qquad\text{and}\qquad
  a_3''X_3+\frac{1}{2}\sum_{i=1}^2\sum_{j=1}^2b_{i,j}''X_iX_j,
  \qquad\text{respectively.}
 $$
 Then the convolution \ $\mu'*\mu''$ \ is a Gaussian measure on \ $\HH$ \
 if and only if there exist \ $t',t''\geq0,$ \ a Gaussian semigroup
 \ $(\mu_t)_{t\sgeq 0}$ \ with infinitesimal generator
 \eqref{STABIL} and an element \ $x\in\HH$ \ which is contained in the center of \
 $\HH$ \ such that either \ $\mu'=\mu_{t'},$ \ $\mu''=\mu_{t''}*\varepsilon_x$ \
 or \ $\mu'=\mu_{t'}*\varepsilon_x,$ \ $\mu''=\mu_{t''}$ \ holds.
 Moreover, in this case \ $a_3=a_3'+a_3''$ \ and \ $b_{i,j}=b_{i,j}'+b_{i,j}'',$
 $1\leq i,j\leq 2.$}

The proof of this statement can be carried out in a direct way
applying Theorem 7.3 in Pap \cite{PAP02}, and Lemma \ref{GAUSSPAR2}
and Proposition \ref{CENTRALGAUSS} of the present paper.
\end{Rem}

\vskip0.1cm

\noindent {\bf \Large Acknowledgements}

The authors are very grateful to the unknown referees for their
important observations and remarks that helped to improve the
presentation of the paper.

The authors have been supported by the Hungarian Scientific Research
Fund under Grant No.~OTKA--T048544/2005. The first author has been
also supported by the Hungarian Scientific Research Fund under Grant
No.~OTKA--F046061/2004.

\parbox{8cm}{M\'aty\'as Barczy\\
             Faculty of Informatics\\
             University of Debrecen\\
             Pf.12\\
             H--4010 Debrecen\\
             Hungary\\ \\
             barczy@inf.unideb.hu}
\hfill
\parbox{5cm}{Gyula Pap\\
             Faculty of Informatics\\
             University of Debrecen\\
             Pf.12\\
             H--4010 Debrecen\\
             Hungary\\ \\
             papgy@inf.unideb.hu}

\end{document}